\documentclass[11pt]{article}
\usepackage[english]{babel}
\usepackage{amsmath, amsfonts, amssymb, amsthm}
\usepackage{graphicx}
\usepackage{array}

\newcommand \R {\mathbb{R}}
\newcommand \N {\mathbb{N}}

\newcommand \la {\langle}
\newcommand \ra {\rangle}
\newcommand \ind {\mathbf{1}}
\newcommand \E {\mathbb{E}}
\def\P{\mathbb{P}}
\newcommand \e {\varepsilon}
\newcommand \und {\underline}
\newcommand \cO {\mathcal{O}}
\def\cA{\mathcal A}
\def\cM{\mathcal M}
\def\w2{w^{(2)}}
\def\wa{w^{(\alpha)}}
\def\Xa{X^{\alpha}}
\def\pa{p^{\alpha}}
\def\cG{\mathcal G}
\def\cGa{\mathcal{G}^{\alpha}}

\newtheorem{theorem}{Theorem}[section]

\newtheorem{lemma}[theorem]{Lemma}

\newtheorem{prop}[theorem]{Proposition}
\theoremstyle{definition}

\newtheorem{rmk}[theorem]{Remark}

\title{Large scale behaviour of the spatial $\Lambda$-Fleming-Viot process}
\author{N. Berestycki\thanks{NB supported in part by EPSRC grants EP/G055068/1 and EP/I03372X/1}\\Statistical Laboratory\\University of Cambridge\\DPMMS\\Wilberforce Road\\Cambridge CB3 0WB, UK\\~\\
A.M. Etheridge\thanks{AME supported in part by EPSRC Grants EP/E065945/1 and EP/G052026/1}\\
Department of Statistics\\University of Oxford\\1 South Parks Road\\Oxford OX1 3TG, UK\\~\\
A. V\'eber\thanks{AV supported in part by the {\em chaire Mod\'elisation Math\'ematique et Biodiversit\'e} of Veolia
Environnement-\'Ecole Polytechnique-Museum National d'Histoire Naturelle-Fondation X and by the
ANR project MANEGE (ANR-090BLAN-0215).}\\
CMAP -\'Ecole Polytechnique\\Route de Saclay\\91128 Palaiseau Cedex\\France}

\begin{document}
\maketitle
\newpage
\begin{abstract}
We consider the spatial $\Lambda$-Fleming-Viot process model (\cite{BEV2010}) for
frequencies of genetic types in
a population living in $\R^d$, in the special case in which there are just two types of individuals,
labelled $0$ and $1$.  At time zero, everyone in a given half-space
has type $1$, whereas everyone in the complementary half-space has
type $0$.  We are concerned with patterns of frequencies of the two types at large space and
time scales.  We consider two cases, one in which the dynamics of the process are driven by
purely `local' events and one incorporating large-scale extinction recolonisation events.  We choose
the frequency of these events in such a way that, under a suitable rescaling of space and time,
the ancestry of a single individual in the
population converges to a symmetric stable process of
index $\alpha\in(1,2]$ (with $\alpha=2$ corresponding to Brownian motion).  We consider the
behaviour of the process of allele frequencies under the same space and time rescaling.  For $\alpha=2$,
and $d\geq 2$ it converges to a deterministic limit.
In all other cases the limit is random and we identify it as the indicator function of a random set.
In particular, there is no local coexistence of types in the limit.  We characterise the set in terms
of a dual process of coalescing symmetric stable processes, which is of interest in its own right.
The complex geometry of the random set is illustrated through simulations.

\bigskip
On \'etudie le processus $\Lambda$-Fleming-Viot spatial (\cite{BEV2010}) mod\'elisant les fr\'equences locales de types g\'en\'etiques dans une population \'evoluant dans $\R^d$. On consid\`ere le cas particulier o\`u il n'y a que deux types possibles, not\'es $0$ et $1$. Initialement, tous les individus pr\'esents dans le demi-espace des points dont la premi\`ere coordonn\'ee est n\'egative sont de type $1$, tandis que les individus pr\'esents dans le demi-espace compl\'ementaire sont de type $0$. On s'int\'eresse au comportement des fr\'equences locales sur des \'echelles de temps et d'espace tr\`es grandes. On consid\`ere deux cas~: dans le premier, l'\'evolution du processus est due uniquement \`a des \'ev\'enements `locaux'; dans le second, on incorpore des \'ev\'enements d'extinction et recolonisation de grande ampleur. On choisit la fr\'equence de ces \'ev\'enements de sorte qu'apr\`es une renormalisation spatiale et temporelle appropri\'ee, la lign\'ee ancestrale d'un individu de la population converge vers un processus $\alpha$-stable sym\'etrique, d'indice $\alpha \in (1,2]$ (o\`u $\alpha=2$ correspond au mouvement brownien). On \'etudie l'\'evolution du processus des fr\'equences all\'eliques aux m\^emes \'echelles spatio-temporelles. Lorsque $\alpha=2$ et $d\geq 2$, celui-ci converge vers un processus d\'eterministe. Dans tous les autres cas, le processus limite est al\'eatoire et on l'identifie comme la fonction indicatrice d'un ensemble al\'eatoire \'evoluant au cours du temps. En particulier, les deux types ne coexistent pas \`a la limite. On caract\'erise chaque ensemble en termes d'un processus dual constitu\'e de mouvements stables sym\'etriques coalescents ayant un int\'er\^et en eux-m\^emes. La g\'eom\'etrie complexe des ensembles limites est illustr\'ee par des simulations.
\end{abstract}



\newpage

\section{Introduction}\label{section: intro}

In this article, we are interested in the behaviour over large space and
time scales
of the \emph{spatial $\Lambda$-Fleming-Viot process} (or SLFV) on $\R^d$.
This process arises as a particular instance of the framework introduced in
\cite{etheridge:2008, BEV2010, barton/kelleher/etheridge:2010}
for modelling allele frequencies (that is frequencies of different
genetic types) in a population that evolves in a spatial continuum.
From the modelling perspective, this framework is interesting as it
overcomes an obstruction to modelling biological populations in continua, dubbed `the pain in the torus' by
Felsenstein (\cite{felsenstein:1975}), which is typified by
the `clumping and extinction' seen in spatial branching process models in
low dimensions.   The key idea of the SLFV framework
is to base reproduction events on a
space-time Poisson process rather than on individuals in the population.
In this way one can define what can be thought of as a continuum version of
the Kimura stepping stone model (\cite{KIM1953}) which is a widely accepted model for
evolution of allele frequencies in spatially {\em subdivided} populations.
Moreover, one can incorporate
large-scale extinction-recolonisation events through a
series of `local' population bottlenecks, each affecting substantial portions of the species range.
Such events dominate the demographic history of many species and, as
we shall see in our results here, can
have a very significant influence on patterns of allele frequencies.

From a mathematical perspective, the SLFV process is a natural extension
to the spatial context of the generalised Fleming-Viot processes which can be traced
to \S3.1.4 of \cite{donnelly/kurtz:1999} but were first studied in detail by
Bertoin \& Le Gall (\cite{bertoin/legall:2003}).  These processes are dual
to the so-called $\Lambda$-coalescents which have been the subject of
intensive study since their introduction over a decade ago
by Donnelly \& Kurtz, Pitman
and Sagitov (\cite{donnelly/kurtz:1999, pitman:1999, sagitov:1999}).
The duality with the generalised Fleming-Viot processes
extends that between the Kingman coalescent and
the Wright-Fisher diffusion and our work here will exploit a similar
duality between
spatial versions of the $\Lambda$-coalescents and the SLFV.  One of the attractions
of the SLFV processes is that they allow us to capture many of the
features of Wright-Fisher noise, but in {\em any} spatial dimension (whereas
stochastic partial differential equations driven by Wright-Fisher noise only
make sense in one dimension).  Thus, although they were originally motivated
by purely biological considerations,
we believe that these models are also of intrinsic mathematical interest.

\subsection{The spatial $\Lambda$-Fleming-Viot process}\label{subs: model}

First we describe the model.
Each individual in the population is assigned a genetic type, from
a compact space $K$, and a location, in $\R^d$.  At time $t$, the population
is represented by a measurable function
$\rho_t:\R^d \rightarrow \cM_1(K)$, where $\cM_1(K)$ is the set of all
probability measures on $K$.  (In fact, as explained in \S\ref{section: duality}, in defining the state space,
$\Xi$, of the process we identify any two such functions that are equal for
Lebesgue-a.e.~$x\in\R^d$.)
The interpretation of the model is as follows: the
population density is uniform across $\R^d$ and, for each $x\in\R^d$, if we
sample an individual from $x$, then its genetic type is determined by
sampling from the probability measure $\rho_t(x)$.

The dynamics of the population are driven by a Poisson point process, $\Pi$, on
$\R\times\R^d\times (0,\infty)\times [0,1]$,
each point of which specifies a (local) extinction-recolonisation event.
If $(t,x,r,u)\in \Pi$, then, at time $t$:
\begin{enumerate}
\item An extinction-recolonisation event affects the closed ball
$B(x,r)\subseteq \R^d$, and nothing happens outside this region.
\item A parent is chosen uniformly in the ball; that is,
we sample a location $z$ uniformly at random over $B(x,r)$ and a type
$k$ according to the distribution $\rho_{t-}(z)$.
\item  For each $y\in B(x,r)$ (including $z$), a fraction $u$ of the local
population is replaced by offspring, whose type is that of the chosen parent.
That is,
$$
\rho_t(y):= (1-u)\rho_{t-}(y) + u\, \delta_k.
$$
\end{enumerate}
Here, we are thinking of reproduction events as
equivalent to (frequent) small-scale extinction-recolonisation events.

In \cite{BEV2010}, the intensity measure of the
Poisson point process $\Pi$ has the form $dt\otimes dx\otimes \zeta(dr,du)$,
thus allowing the `impact', $u$, of a reproduction event to depend on its
radius, $r$.  For instance, small-scale reproduction events may affect only
a tiny fraction of individuals, compared to massive extinction-recolonisation
events which could wipe out most of the population in a large
geographical region.
Of course, we require some conditions on the intensity of $\Pi$ if our
process is to be well-defined: according to Theorem~4.2 of \cite{BEV2010}
(stated for $d=2$, but the proof is identical for any dimension $d\geq 1$),
the corresponding spatial $\Lambda$-Fleming-Viot process is well-defined
whenever
\begin{equation}\label{cond existence}
\int_{(0,\infty)\times [0,1]}
\zeta(dr,du) \; uV_r <\infty,
\end{equation}
where $V_r$ denotes the volume of a $d$-dimensional ball of radius $r$.

\subsection{Main results}\label{subs: frameworks}

Our previous mathematical analysis of the SLFV process (\cite{BEV2010, etheridge/veber:2011})
has been concerned with understanding the genealogical relationships
between individuals sampled from the population.  Here,
although studying the lineages ancestral to a sample from the population
will be fundamental to our analysis,  we are interested
in understanding the
patterns of allele frequencies that result from such a model.

Saadi (\cite{saadi:2011}) considers a closely related model (which differs
from ours only in that the location of the `parent' in a reproduction event
is always taken to be the centre of the event).  He considers the most
biologically interesting case of two
spatial dimensions and, for simplicity, takes all reproduction events to
have fixed size $r$ and fixed impact $u\in(0,1]$.  Notice that if
a particular genetic type is present in a region at some time $t$, then,
unless $u=1$, it will also be there at all later times.
Saadi shows that if
a particular genetic type is only present in a
bounded region at time zero, then, with probability one, its {\em range},
that is the region in which it is {\em ever} seen is bounded.
On the other hand, the {\em shape} of this region will be complex.  In
order to try to gain some understanding of the boundary of the range,
he has also simulated a simpler situation.
The idea is to consider just two `competing' types on a two-dimensional
torus which we can identify with $(-L,L]^2\subseteq\R^2$.
At time zero, all points of the torus with a non-positive first
coordinate are of one type and all with a strictly positive first coordinate
are of the other type.
The region in which both types coexist develops in a rather
complicated way, but it is natural to ask whether if one `stands back' and
views the process over large spatial scales (at sufficiently large
times) a simpler pattern emerges.  Saadi's simulations were the starting
point for our work here.

We shall concentrate our attention on two special cases of the SLFV model,
in both of which individuals can be one of only two genetic types, labelled $0$ and $1$.
Evidently it is then enough to consider the proportion of type\,-1 individuals at each site and so
we define, for every $x\in \R^d$ and $t\geq 0$,
\begin{equation}\label{def density}
w(t,x):= \rho_t(x)(\{1\}).
\end{equation}
For the sake of clarity, we shall also take the fraction $u\in (0,1]$ to
be the same for all events. In our previous notation, this corresponds to
taking $\zeta(dr,dv)=\mu(dr)\delta_u(dv)$, for a measure $\mu$ on
$(0,\infty)$.  We shall allow the measure $\mu$ to take two forms:
\begin{description}
\item[Case A (fixed radius):] We fix $r\in (0,\infty)$, and choose $\mu$ to be the Dirac mass at $r$.
\item[Case B (heavy-tailed distribution):] We fix $\alpha \in (1,2)$ and define the measure $\mu$ by
\begin{equation}\label{def mu B}
\mu(dr) = r^{-\alpha-d-1}\ind_{\{r\geq 1\}}\, dr,
\end{equation}
where we recall that $d$ is the dimension of the geographical space.
\end{description}
It is easy to check that the condition~(\ref{cond existence}) which
guarantees existence of the SLFV process is satisfied in both cases.

Case A bears some similarity to the nearest-neighbour voter model,
in that an individual spreads its type (/opinion) in a `close' neighbourhood.
Case B incorporates some large-scale events and consequently, as we shall
see, behaves
very differently.  The particular form of $\mu$ is motivated by the fact
that with this choice, under a suitable rescaling of space and time,
the motion of an ancestral lineage will converge to a symmetric
$\alpha$-stable L\'evy process (and, more generally, the ancestry of finitely
many individuals converges to a system of coalescing dependent $\alpha$-stable
processes, see \S\ref{S:heavy}).  Combined with duality, this
will imply that with the
same space-time rescaling, the forwards in time process of allele frequencies
will also converge to a non-trivial limit.

Suppose that the initial condition of the process is
$$
w(0,x) = \ind_{\{x_{(1)} \leq 0\}},
$$
where here again $x_{(1)}$ denotes the first coordinate of $x$.  In words,
we start from a half-space $H$ of 1's.  Let us set $\alpha=2$ in Case A, and, for a given $\alpha\in (1,2]$
and any $n\in \N$, define the rescaled
density $w^n$ by
$$
w^n(t,x):= w(nt,n^{1/\alpha}x),\qquad \qquad t\geq 0,\ x\in \R^d.
$$
We denote by $\rho^n$ the $\Xi$-valued process whose local
density of $1$'s at time $t$ is $w^n(t,\cdot)$. Our main results are the
following two theorems, which describe the asymptotic behaviour of
$\rho^n$ as $n$ tends to infinity.
In Case A, $\sigma^2$ is the variance of the displacement, after one unit of time, of a single ancestral lineage
from its starting point
(see~(\ref{def sigma})).

\begin{theorem}\label{theo: A}{\bf (Case A)}
There exists a $\Xi$-valued process $\{\rho^{(2)}_t,\, t\geq 0\}$ such that
$$
\rho^n \longrightarrow \rho^{(2)} \qquad \mathrm{as}\ n\rightarrow \infty,
$$
in the sense of weak convergence of the (temporal) finite-dimensional distributions.

Furthermore, at every time $t\geq0$, the local density $w^{(2)}(t,\cdot):=\rho^{(2)}_t(\{1\})$ of type\,-$1$
individuals can be described as follows. If $X$ denotes standard $d$-dimensional Brownian motion and
$$
p^2(t,x):= \P_x[X_{\sigma^2t}\in H],\qquad t\geq 0,\ x\in\R^d,
$$
then:
\begin{enumerate}
\item{If $d=1$, for every $t\geq 0$ and a.e. $x\in\R$, $\w2(t,x)$ is a Bernoulli random variable
with parameter $p^2(t,x)$.  The correlations between their values at distinct sites of $\R$ are non-trivial
and are described in (\ref{correl A1}).}
\item{If $d\geq 2$, for every $t\geq 0$ and a.e. $x\in \R^d$, $\w2(t,x)$ is deterministic and equal to $p^2(t,x)$.}
\end{enumerate}
\end{theorem}

\begin{rmk}
Note that, in one dimension, the two types almost surely do not coexist
at any given point, since $\w2(t,x)$ is a Bernoulli random variable.
However, in higher dimensions, the two types 0 and 1 do coexist at every
site instantaneously.
\end{rmk}

\begin{rmk}
Although we have expressed everything in terms of densities, the convergence
in Theorem~\ref{theo: A}, which we define explicitly in \S\ref{section: duality},
is equivalent the convergence for all $j \geq 1$ of the measures
$$
\rho_t^n(x_1, d \kappa_1) \ldots \rho_t^n (x_j, d \kappa_j) dx_1 \ldots dx_j.
$$
(This is the vague convergence of nonnegative Radon measures on $(\R^d \times \{0,1\})^j$.)
See \cite{VW2011} for a measure-valued formulation of the SLFV and for a
proof of this equivalence.
\end{rmk}

\begin{rmk}
The quantity $p^2(t,x)$ implicitly depends on the dimension. Also, since $u$ and $r$
are fixed, substituting in~(\ref{def sigma}),
$$
\sigma^2 = \frac{u}{d V_r}\int_{\R^d}dz\ |z|^2  L_r(z) \qquad \left(= \frac{4u r^3}{3}\ \mathrm{when}\ d=1\right)
$$
is finite and proportional to $u$.  Indeed, $L_r(z):= \mathrm{Vol}(B(0,r)\cap B(0,z)) = (2r-|z|)_+$ in
dimension $1$ and, more generally, $L_r(z)\leq \ind_{\{|z|\leq 2r\}}V_r$
for any $d\geq 1$.
\end{rmk}

In contrast to the case of fixed radii, in Case B, in the limit as
$n\rightarrow\infty$ types are always
segregated, irrespective of dimension.

\begin{theorem}\label{theo: B}{\bf (Case B)}
There exists a $\Xi$-valued process $\{\rho^{(\alpha)}_t,\, t\geq 0\}$ such that
$$
\rho^n \longrightarrow \rho^{(\alpha)}\qquad \mathrm{as}\ n\rightarrow \infty,
$$
in the sense of weak convergence of the (temporal) finite-dimensional distributions.

Furthermore, there exists a symmetric $\alpha$-stable process $\Xa$ such that if
$$
\pa(x,t):= \P_x\big[\Xa_{ut}\in H\big],\qquad t\geq 0,\ x\in \R^d,
$$
then for every $t> 0$ and a.e. $x\in \R^d$, $\wa(t,x)$ is a Bernoulli random variable
with parameter $\pa(t,x)$.
The correlations between the values of the densities at different sites (and at the same time $t$) are again
given by (\ref{correl A1}) (or (\ref{point duality})), where the process $\xi^{\infty}$ is now the system of
coalescing $\alpha$-stable processes obtained in Proposition~\ref{prop: geneal B}.
\end{theorem}

Here again, one should notice that the speed of evolution of the limiting
process is proportional to the parameter $u$.

Comparing the results of Theorem~\ref{theo: A} and Theorem~\ref{theo: B}, one
can see that very large extinction-recolonisation events create
correlations between local genetic diversities over a much larger spatial
scale ($n^{1/\alpha}\gg \sqrt{n}$) than purely local reproduction events.
This is because an ancestral lineage can move a
distance $\cO(n^{1/\alpha})$ over the course of
$n$ generations.
One might initially
guess that, since the motion of a single ancestral lineage under our
rescaling converges to a symmetric stable process,
two distinct ancestral lineages would (asymptotically) only meet (and
thus have a chance to coalesce) in dimensions where the stable process
hits points.  This is precisely what we see in Case A and, in that case,
lies behind the
deterministic limit in $d\geq 2$.  However, this
is where the dependence between ancestral lineages in the SLFV process
(see \S\ref{subs: duality})
comes into play.  The
detailed analysis of the ancestral process for Case B (which we
present in \S\ref{S:heavy})
reveals that `very large' events are frequent enough to capture
lineages that have moved to arbitrarily large separations.
In particular, Lemma~\ref{L:coalH} shows
that, in Case B, any finite sample of individuals will find its most
recent common ancestor in finite time a.s. (see also Remark~\ref{rk: MRCA}).
The large events will, momentarily, create extensive areas in which
the two genetic types coexist.  Our analysis will also show that, under our
rescaling, `small' events then
occur sufficiently quickly to instantaneously restore the allele
frequencies in each infinitesimal region
to $0$ or $1$ (see also
the simulations presented in \S\ref{simulations}).

The rest of the paper is laid out as follows.  In \S\ref{simulations},
we present some simulations that illustrate the results and the mechanisms underlying them.
In \S\ref{section: duality},
we are explicit about the meaning of `weak convergence of the
(temporal) finite-dimensional distributions' and we describe the duality
between allele frequencies and
ancestral processes that provides the main tool in our proofs.
It is then used to find conditions, expressed in terms of the genealogical trees
relating individuals in a sample from the population, under which $\wa(t,x)$ (at each time $t> 0$ and
a.e. point  $x\in\R^d$) takes the particular forms seen in our main theorems
(see Lemma \ref{lem: density}).
Theorems~\ref{theo: A} and \ref{theo: B} are then proved in
\S\ref{section: proof A} and \S\ref{S:heavy} respectively.  This last
section also contains some results
(Lemma~\ref{L:coalH} and the accompanying remark),
of independent interest, on the system of coalescing (dependent) L\'evy processes that
generates the genealogical trees relating
a sample of individuals from the limiting population.

\section{Simulations}
\label{simulations}

Our results show that in the cases where the
rescaled density of type $1$ individuals converges to a {\em random} limit,
at any fixed time that limit takes the form of the indicator function of a random set.  In one dimension, provided that
either $u=1$ or $\alpha=2$ (the radius of events is fixed), the set takes a simple form, but for $\alpha\in (1,2)$ this is no
longer the case.
In this section we present some simulations that illustrate the complex geometry of the limiting random sets and
the mechanism that leads to their creation.
We are extremely grateful to Jerome Kelleher from the University of
Edinburgh for performing these simulations and producing the figures.

First suppose that we are in one spatial dimension.
If $u=1$, then at every stage of the rescaling
we will have $w^n(t,x)=\mathbf 1_{I^n_t}(x)$ where $I^n_t$ is a half-line with right endpoint $R_t$ following a
random walk on $\R$.  Under our rescaling, as $n\rightarrow\infty$, the process $R_t$ will converge to a Brownian
motion if $\alpha=2$ and to a symmetric stable process of index $\alpha$ for $\alpha\in (1,2)$.
If $\alpha=2$, and $d=1$, then the same is true for $u<1$.  This can be understood via the dual
process of ancestral lineages.  As we shall see, this converges to a system of independent Brownian motions which coalesce
instantaneously on meeting.  The type of an individual sampled at $x$ at time $t$ is determined by the type at time
$t$ before the present of
the corresponding ancestral lineage.  Since the Brownian motions are continuous, and they coalesce as soon as they meet,
it is impossible
for two lineages to `cross over'.  Consequently, asymptotically, if a lineage started from $x$ traces
back to a point to the left of the origin at time $t$ before the present, then so must all lineages
started from points to the left of $x$.  As a result, at time $t$, the density of type $1$s will still
be the indicator function of a half-line.  The boundary, $R_t$, moves in the same
way as a single ancestral lineage, that is as a Brownian motion with a clock that runs at a rate proportional to $u$.
Figure~\ref{fig: d=1, a=2} shows the results of a simulation of the process of allelic types in this case.
In two dimensions, two Brownian motions won't meet and so for $\alpha=2$, asymptotically,
the ancestral lineages will just look like
independent Brownian motions and forwards in time, asymptotically, allele frequencies are smeared out by
the deterministic heat flow.
\begin{figure}[htp]
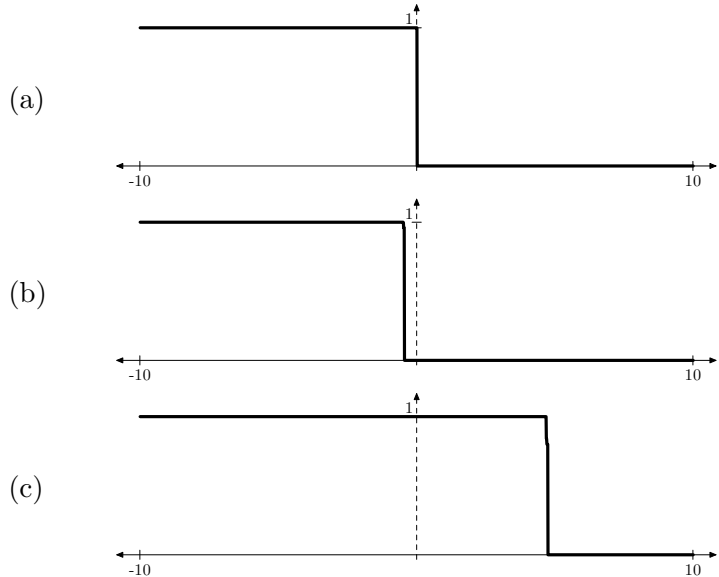

\centering
\newcommand{\figwidth}{0.5\textwidth}
\begin{tabular}{m{1cm}m{0.5\textwidth}}
(a) & \includegraphics[width=\figwidth]{d1fr.0} \\
(b) & \includegraphics[width=\figwidth]{d1fr.1} \\
(c) & \includegraphics[width=\figwidth]{d1fr.2} \\
\end{tabular}
\caption{Fixed radius in $d=1$ on a line of length 20. (a) initial conditions; (b)
after $10^5$ events;
(c) after $10^7$ events. The model parameters are $u = 0.8$, $r = 0.033$, $n=10^3$.
}
\label{fig: d=1, a=2}
\end{figure}

The case $\alpha\in(1,2)$ is much more interesting.  Now, even in the limit, ancestral lineages evolve in a series
of jumps and if $u<1$ they {\em can} `cross over'.  Thus although our results show that the limiting allele frequencies always
look like the indicator function of a random set, even in $d=1$ we can no longer expect that set to be
connected.  Forwards in time what our results suggest, and simulations confirm, is that a large
event can create a region in which allele frequencies are strictly between zero and one, but these
frequencies are rapidly (and asymptotically instantaneously) `resolved' by `small' events so that the
state is restored to being the indicator function of a set.
Figure~\ref{fig: d=1, a<1} shows how on the line this mechanism leads to allele frequencies that look like a series of `crenellations'.
Even in one spatial dimension, our methods
are not powerful enough to allow us to capture detailed information about the random sets observed
in the limit.
\begin{figure}[htp]
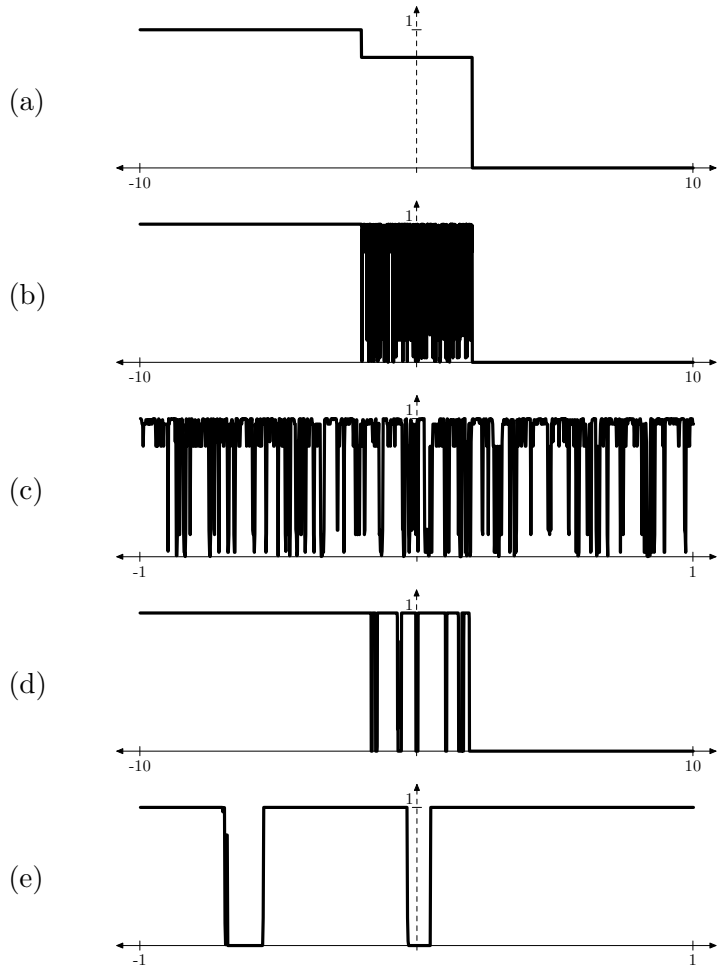

\centering
\newcommand{\figwidth}{0.5\textwidth}
\begin{tabular}{m{1cm}m{0.5\textwidth}}
(a) & \includegraphics[width=\figwidth]{d1vr.0} \\
(b) & \includegraphics[width=\figwidth]{d1vr.1} \\
(c) & \includegraphics[width=\figwidth]{d1vr.2} \\
(d) & \includegraphics[width=\figwidth]{d1vr.3} \\
(e) & \includegraphics[width=\figwidth]{d1vr.4} \\
\end{tabular}
\caption{Variable radius in $d=1$ on a line of length 20. (a) initial conditions; (b)
after $100$ events, full range; (c) after $100$ events, zooming in;
(d) after $10^6$ events, full range;
(e) after $10^6$ events, zooming in. The model parameters are $u = 0.8$, $n = 10^4$ and
$\alpha = 1.3$.
}
\label{fig: d=1, a<1}
\end{figure}

Figure~\ref{fig: d=2, a<1} illustrates the same mechanism in two spatial dimensions.
To isolate the effect in which we are interested, we
suppose that a large event covers a previously unblemished portion of the interface and
observe the resolution of the resulting patch of coexistence.

\begin{figure}[htp]
\centering
\begin{tabular}{@{}c@{}c@{}c@{}}
\includegraphics[width=0.33\textwidth]{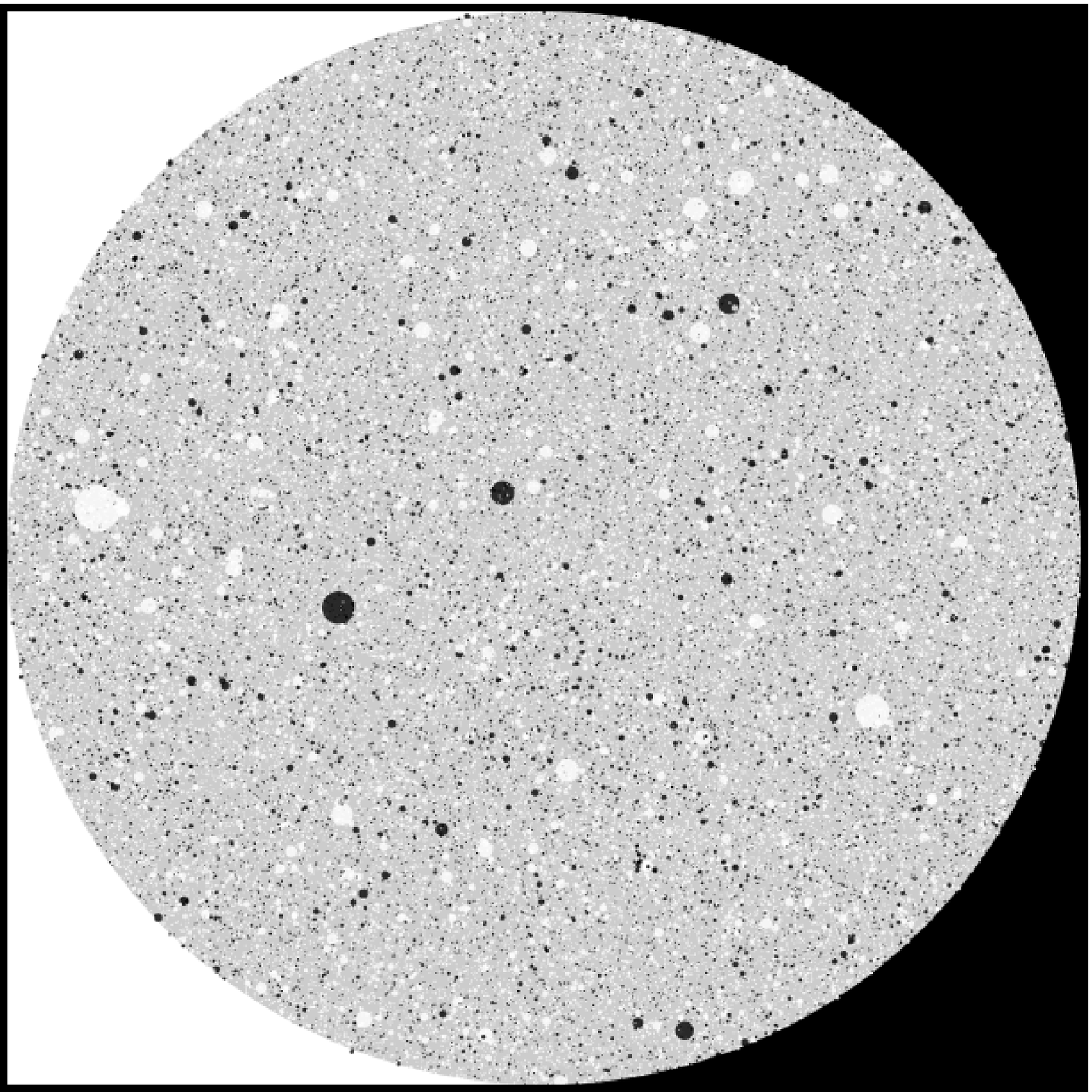} &
\includegraphics[width=0.33\textwidth]{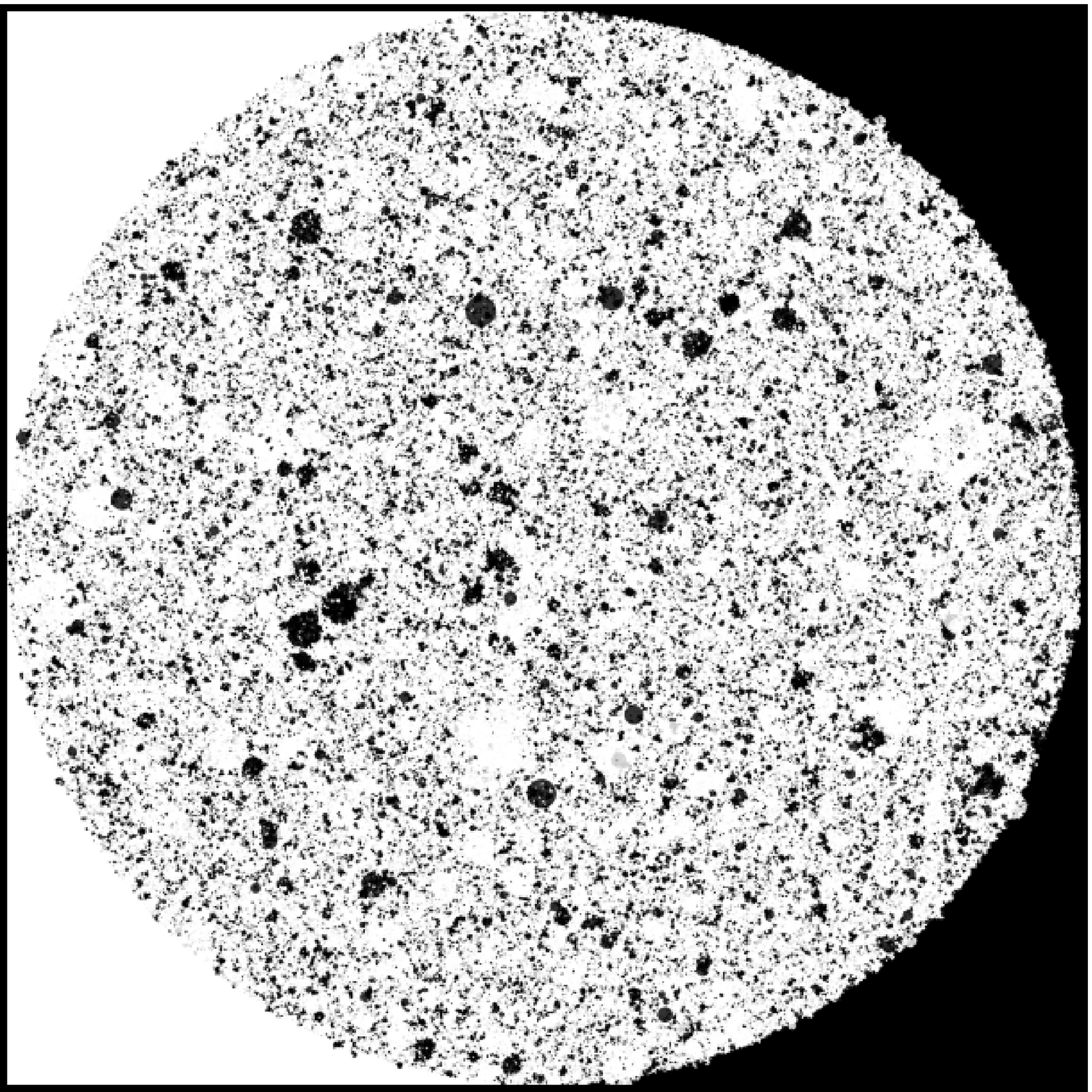} &
\includegraphics[width=0.33\textwidth]{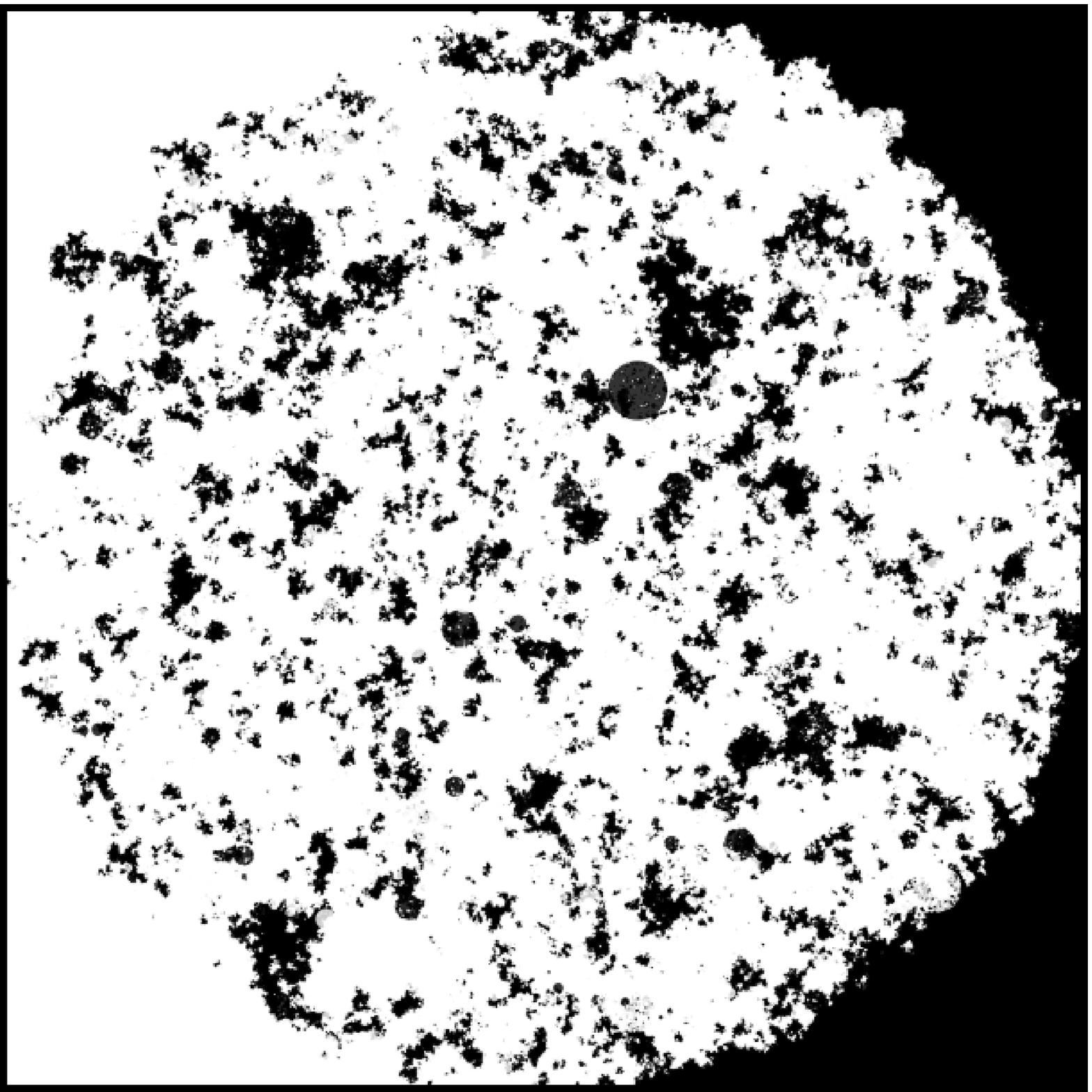} \\
(a) & (b) & (c) \\
\end{tabular}
\caption{Model in $d=2$ after (a) $10^5$; (b) $10^6$; and (c) $10^7$ events.
We have a square range of edge $8$, and the initial patch is a circle
of radius 4 with frequency $0.8$ (white is frequency 1, black is 0).
The model parameters are $u = 0.8$, $\alpha=1.3$ and $n = 10^3$.
}
\label{fig: d=2, a<1}
\end{figure}

\section{Convergence and duality}\label{section: duality}

\subsection{State-space and form of convergence}\label{subs: form conv}

In order to make the convergence in Theorems~\ref{theo: A} and \ref{theo: B} explicit, let us recall
some facts about the state space of the SLFV from \cite{BEV2010}. In \S\ref{subs: model}, we described the process
as taking its values in the set $\tilde{\Xi}$ of all measurable functions $\rho:\R^d \rightarrow \cM_1(K)$
(where the compact type space $K$ is now $\{0,1\}$). In fact, we need to define an equivalence
relation on this space by setting
$$
\rho \sim \rho'\qquad \Leftrightarrow \qquad \mathrm{Vol}\big(\{x\in \R^d\,:\, \rho(x)\neq \rho'(x)\}\big)=0.
$$
The state-space $\Xi$ of the SLFV is then defined as the quotient space $\tilde{\Xi}/\sim$ of equivalence classes
of $\sim$.

$\Xi$ is endowed with a natural topology described in $\S 3$ of \cite{EVA1997}. As mentioned earlier, it is possible to check that this topology is equivalent to the topology of vague convergence of the measures $\rho(x_1, d \kappa_1) \ldots \rho(x_j, d\kappa_j) dx_1 \ldots dx_j$ for all $j\geq 1$ as Radon measures on $(\R^d \times \{0,1\})^j$.
In practice, Lemma~4.1 in \cite{BEV2010} provides us with a family of functions which is dense in the set of continuous functions
on $\Xi$.  To introduce this class of functions, for a space $E$, let $C(E)$ denote the set of all continuous functions on
$E$ and, for a measure $\nu$, let $L^1(\nu)$ be the set of all functions which are integrable with respect to $\nu$.
For every $j\in \N$, $\psi\in C((\R^d)^j)\cap L^1(dx^{\otimes j})$ and $\chi_1,\ldots,\chi_j\in C(K)$,
we define the function $I_j(\cdot\,;\,\psi,(\chi_i)_{1\leq i\leq j})$ as follows.  For every $\rho\in \Xi$,
$$
I_j(\rho\,;\,\psi,(\chi_i)_{1\leq i\leq j}):= \int_{(\R^d)^j}dx_1\ldots dx_j\
\psi(x_1,\ldots,x_j)\left(\prod_{i=1}^j\big\la \chi_i, \rho(x_i)\big\ra\right),
$$
where $\la f,\nu\ra$ is the integral of the function $f$ with respect to the measure $\nu$.

Since in our setting $K=\{0,1\}$, we have, for every $\chi$,
\begin{align*}
\la \chi,\rho(x)\ra & =
\chi(1) w(x) + \chi(0)(1-w(x))\\
&= \big(\chi(1)- \chi(0)\big) w(x) +\chi(0),
\end{align*}
where, as before, $w(x):=\rho(x)(\{1\})$ denotes the mass of $1$'s at site $x$.
We can therefore restrict our attention to the set of functions $I_j$ such that
$\chi_i = \mathrm{Id}$ for every $i\in \{1,\ldots,j\}$, that is
\begin{equation}\label{test functions}
I_j(\rho\,;\, \psi)= \int_{(\R^d)^j}dx_1\ldots dx_j\
\psi(x_1,\ldots,x_j)\left(\prod_{i=1}^j w(x_i)\right).
\end{equation}
Indeed, any $I_j(\cdot\,;\,\psi,(\chi_i)_{1\leq i\leq j})$ can be written as a finite linear
combination of functions of the
form~(\ref{test functions}). The convergence stated in Theorems~\ref{theo: A} and \ref{theo: B} can now be expressed
for a single time $t\geq 0$ as: for every $j\in \N$ and $\psi\in C((\R^d)^j)\cap L^1(dx^{\otimes j})$,
\begin{equation}\label{sense conv}
\lim_{n\rightarrow \infty}\E\big[I_j\big(\rho^n_t\,;\, \psi\big)\big] = \E\big[I_j \big(\rho^{(\alpha)}_t\,;\, \psi\big)\big].
\end{equation}
The extension of this definition of convergence to joint convergence at several times $t_1,\ldots,t_k$ is straightforward.

\subsection{Duality between the SLFV and its genealogies}\label{subs: duality}

The proofs of Theorems~\ref{theo: A} and \ref{theo: B} rely on
a duality relation between the SLFV process, and the system of coalescing jump processes that we call the
\emph{genealogical process of a sample of individuals} from the population.  We recall this relation in the particular form in which we
shall need it.  In particular we restrict our attention to $K=\{0,1\}$.  A more general statement
(and proofs) can be found in \S4 of \cite{BEV2010}.

First suppose that we wish to trace
the ancestry of a single individual alive in the current population.
Let us, for now, work in a general setting as it will allow us to
understand condition~(\ref{cond existence})
a little better.
Since the
model is translation invariant, without loss of generality we may suppose that
the individual is currently at the origin in $\R^d$.  Tracing backwards
in time, at the first time in the past when $0$ is in the area $B(x,r)$
affected by a reproduction event, our individual has probability $u$ of
being an offspring of that event, in which case the ancestral lineage jumps
to the position of the parent (which is uniformly distributed
on $B(x,r)$).  Since the Poisson process driving events is reversible, we
see that the rate at which our ancestral lineage experiences a jump is
$$
\int_{\R^d}\int_{(0,\infty)\times [0,1]}\zeta(dr,du)dx \;
\ind_{\{0\in B(x,r)\}}u
= \int_{(0,\infty)\times [0,1]} \zeta(dr,du)\; uV_r.
$$
By translation invariance in time and space of the law of $\Pi$, this tells us that the
quantity in (\ref{cond existence}) is just the instantaneous jump rate
of an ancestral lineage (at any time and any location), and we are
requiring it to be finite. We refer to \S4 in \cite{BEV2010} for an
explanation of why this guarantees existence and uniqueness of the
process $(\rho_t)_{t\geq 0}$.

We will need to be more precise about the
law of the compound Poisson process followed by an ancestral lineage and so
we now
establish the rate at which it jumps from $0$ to $z$ (or, by
translation invariance, from $y$ to $y+z$).
In order for such a jump to occur, first $0$ and $z$ must both belong to
the area hit by the event; second our lineage at $0$ must belong to the
fraction $u$ of individuals replaced; and third the parent must be chosen
from site $z$. The intensity measure of the jump process is
therefore equal to \begin{eqnarray}
m(dz) &:= &\int_{\R^d}\int_{(0,\infty)\times [0,1]}\zeta(dr,du)dx
\,\ind_{\{x\in B(0,r)\cap B(z,r)\}}u\,\frac{dz}{V_r}\,
\nonumber \\
&=& \left(\int_{(0,\infty)\times [0,1]} \ \zeta(dr,du)
\frac{u L_r(z)}{V_r}
\right) dz,
\label{jump intensity}
\end{eqnarray}
where $L_r(z)$ denotes the volume of the intersection $B(0,r)\cap B(z,r)$.
(To see this, note that for an event of radius $r$ to affect both $0$ and $z$,
its centre, $x$, must
lie in the region $B(0,r)\cap B(z,r)$ and that since the parent is chosen
uniformly from the region, the factor $1/V_r$ arises as
the density of the uniform distribution on $B(x,r)$.)
In particular, by rotational symmetry, in the special case where the variance of the
displacement of a lineage over one unit of time
is finite,
its covariance matrix is of the
form $\sigma^2 \mathrm{Id}$, with
\begin{equation}\label{def sigma}
\sigma^2:= \int_{\R^d}m(dz)\ \big(z_{(1)}\big)^2 = \frac{1}{d}\int_{\R^d}m(dz)\ |z|^2
\end{equation}
(here, $z_{(1)}$ denotes the first
coordinate of $z$, and $|z|$ its $L^2$-norm).

Much of our analysis will rest upon understanding the ancestry of (larger) samples from
the population, and these can be established in much the same way as the motion
of a single ancestral lineage.
If we sample $k$ individuals (possibly from the same location),
the ancestry is given by a system of (finite-rate) jump processes, which
are \emph{a priori} correlated, since their jumps are generated by the same
Poisson point process of events. Furthermore, if at least two of them are
encompassed by the same event and lie within the fraction of the local
population replaced, then these lineages trace back to the same parent and
thus merge into a single lineage during the event.
Tracing further back in time, that single lineage and all other remaining
lineages continue to evolve in the same manner.
Note that if $u<1$, there may be other lineages in
the ball where the event takes place, but not in the sub-population replaced.
Such lineages neither jump nor coalesce during the event.

Let $(\cA_t)_{t\geq 0}$ be a system of finitely many (the initial number will always be specified explicitly)
ancestral lineages as described above.  That is, each lineage follows a finite-rate jump process
with jump intensity~(\ref{jump intensity}),
and two or more lineages coalesce whenever they are
affected by the same event. See Equation~(\ref{generator A}) in \S\ref{S:heavy} for an expression for the generator
of this process, in the particular case where $u$ is fixed.  For every $t\geq 0$, let us write $N_t$ for the number of
distinct lineages at time $t$, and $\xi_t^1,\ldots,\xi_t^{N_t}\in \R^d$ for their spatial locations at that time.
The weak duality relation we shall use in the sequel is also based on the family of functions in~(\ref{test functions}),
and states that for every $j\in \N$ and $\psi\in C((\R^d)^j)\cap L^1(dx^{\otimes j})$, we have, for every $t\geq 0$, \setlength\arraycolsep{1pt}
\begin{eqnarray}
\int_{(\R^d)^j}&&  dx_1\ldots dx_j \
\psi(x_1,\ldots,x_j)\ \E\big[w(t,x_1)\ldots w(t,x_j)\,\big|\, w(0,\,\cdot)=w_0\big]
\label{eq duality}\\
&=& \int_{(\R^d)^j}
dx_1\ldots dx_j\
\psi(x_1,\ldots,x_j)\ \E\big[w_0(\xi_t^1)\cdots w_0\big(\xi_t^{N_t}\big)\, \big|\, N_0=j,\ \xi_0^1=x_1,\ldots, \xi_0^j=x_j\big].
\nonumber
\end{eqnarray}
Since (\ref{eq duality}) is valid for all functions $\psi$ as above, we also have for Lebesgue-a.e. $(x_1,\ldots,x_j)$,
\begin{equation}\label{point duality}
\E\big[w(t,x_1)\ldots w(t,x_j)\,\big|\, w(0,\,\cdot)=w_0\big]= \E\big[w_0(\xi_t^1)\cdots w_0\big(\xi_t^{N_t}\big)\, \big|\, N_0=j,\ \xi_0^1=x_1,\ldots, \xi_0^j=x_j\big].
\end{equation}
\begin{rmk}
The weak duality in~(\ref{point duality}) is very similar to the duality between the Kimura
stepping stone model and a system of coalescing random walks (see e.g. Chap.6 of \cite{ETH2011}).  Here,
however, in contrast to the discrete space setting, we cannot deduce an expression for the second or higher
order moments of the $w(t,x)$'s since~(\ref{point duality}) only holds for Lebesgue-a.e $j$-tuple $(x_1,\ldots,x_j)$
(and the $x_i$'s are pairwise distinct for Lebesgue-a.e. vector).
The problem stems from the fact the actual object with which we are dealing is the random measure $w(t,x)dx$ and not the collection
$\{w(t,x)\}_{x\in \R^d}$. The topology on $\Xi$ is too weak to consider the evolution of the density of $1$'s at every single point,
and we are obliged to characterize this density through a local averaging procedure, see~(\ref{local average}).
\end{rmk}

Thanks to~(\ref{eq duality}), proving the convergence of $\rho^n_t\equiv \{w(nt,n^{1/\alpha}x)\}_{x\in \R^d}$ boils
down to showing that the genealogical process relating a finite sample of individuals converges, and to transferring
the result to the forwards-in-time process.  In addition, these duality relations enable us to obtain an explicit
description of the local densities $\wa(t,x)$.  Indeed, (\ref{sense conv}) and (\ref{eq duality}) lead us to
an implicit characterisation of the limiting random field $\rho^{(\alpha)}$ through the values of
$$
\E\big[I_j (\rho^{(\alpha)}_t\,;\, \psi)\big] = \E\bigg[\int_{(\R^d)^j}
dx_1\ldots dx_j\
\psi(x_1,\ldots,x_j)\bigg(\prod_{i=1}^j\wa(t,x_j)\bigg)
\bigg].
$$
However, the following result gives us more information on the form of the $\wa(t,x)$'s.
\begin{lemma}\label{lem: density}
Suppose that $(\rho_t)_{t\geq 0}$ is a $\Xi$-valued process dual to an exchangeable and consistent system of coalescing
Markov processes $(\cA_t)_{t\geq 0}$ through the relations (\ref{eq duality}).  Let $(\xi_t)_{t\geq 0}$ denote the
Markov process followed by a single lineage, and suppose that the initial condition of $\rho$ is such that for every $t>0$,
the map $z\mapsto \E_z[w(0,\xi_t)]$ is continuous on $\R^d$ (where as usual $\E_z$ denotes expectation under $\P_z$).
\begin{enumerate}
\item[(i)] If for every $\e>0$ we have
\begin{equation}\label{cond Bernoulli}
\lim_{|y-x|\rightarrow 0}\P\big[\mathrm{lineages}\ 1\ \mathrm{and}\ 2 \ \mathrm{have\ not\ coalesced\ by\ time}\ \e\,\big|\, \xi_0^1 = x, \xi_0^2=y\big] = 0,
\end{equation}
where the convergence is uniform with respect to $x\in\R^d$,
then for every $t> 0$ and a.e. $x\in \R^d$, $w(t,x)$ is a Bernoulli random variable with parameter $\E_x[w(0,\xi_t)]$.
\item[(ii)] If $(\cA_t)_{t\geq 0}$ is a system of independent Markov processes which never coalesce whenever they start from distinct locations, then for every $t>0$ and a.e. $x\in \R^d$, $w(t,x)$ is deterministic and equal to $\E_x[w(0,\xi_t)]$.
\end{enumerate}
\end{lemma}
Here, by `exchangeable' we mean that the law of $(\cA_t)_{t\geq 0}$ is invariant under relabelling of the initial
lineages; `consistent' means that for every $j\in \N$, if $\cA$ starts with $j+1$ lineages but we only follow the
evolution of the first $j$ of them, we obtain a system of coalescing Markov processes that has the same law as
$\cA$ started with only $j$ lineages.  In other words, the evolution of the $(j+1)$-st lineage does not influence
that of the other $j$. It is not difficult to see that the system $(\cA_t)_{t\geq 0}$ introduced at the beginning
of this section is indeed exchangeable and consistent (since each lineage present in the area hit by an event is
affected with probability $u$ independently of all others).  The limiting genealogies we shall obtain will
inherit these properties.

\medskip

\noindent{\bf Proof of Lemma \ref{lem: density}.}
Let us fix $t\geq 0$, and consider the random measure $\ell(dx)$ on $\R^d$ defined by:
for every nonnegative measurable function $\psi$,
\begin{equation}\label{def ell}
\int_{\R^d}
\ell(dx)\ \psi(x)
:= \int_{\R^d}dx\ \psi(x)w(t,x).
\end{equation}
Notice that, according to the description of $\Xi$ given in \S\ref{subs: form conv}, $w(t,\cdot):=\rho_t(\cdot)(\{1\})$
is in fact an equivalence class of functions of the form $\tilde{w}:\R^d\rightarrow [0,1]$.
Two representatives of $w(t,\cdot)$ differ only on a Lebesgue negligible subset of $\R^d$.  For the rest of
this proof we assume that for every $\omega$ in the probability space $(\Omega,\mathcal{F},\P)$ on which $\rho_t$
is defined, we have fixed a representative $\tilde{w}(\omega):\R^d\rightarrow [0,1]$ of $w(\omega,t,\cdot)$ and
define $\ell(\omega,dx)$ as in (\ref{def ell}), with $w(\omega,t,\cdot)$ replaced by $\tilde{w}(\omega,\cdot)$.

Let $(\varphi_m)_{m\in \N}$ be a sequence of continuous functions on $\R^d$ such that for every $m$,
$0\leq \varphi_m\leq 1$, $\varphi_m\equiv 1$ on $B(0,1/m)$ and $\varphi_m\equiv 0$ outside $B(0,2/m)$.
Let us write $\varphi_m(\R^d)$ for the integral $\int_{\R^d}dz\ \varphi_m(z)$.
Since $\tilde{w}$ is locally integrable (it has values in $[0,1]$), the Lebesgue Differentiation Theorem guarantees
that for every $\omega\in \Omega$, there exists a Lebesgue null set $\mathcal{N}(\omega)$ such that for every
$x\notin \mathcal{N}(\omega)$,
\begin{equation}\label{local average}
\lim_{m\rightarrow \infty}\frac{1}{\varphi_m(\R^d)}\int_{\R^d}
\ell(\omega,dz) \
\varphi_m(x+z)
= \tilde{w}(\omega,x).
\end{equation}
Consequently, by Fubini's theorem there exists a Lebesgue null set $\mathcal{O}$ such that for every
$x\notin \mathcal{O}$, the convergence in (\ref{local average}) occurs with $\P(d\omega)$-probability one.
Evidently, if we can show that the random variable $\tilde{w}(x)$ is as in the
statement of Lemma~\ref{lem: density} for every $x\notin \mathcal{O}$, we shall obtain the desired result
for $w(t,\cdot)$.

Now fix $x\in \R^d\setminus\mathcal{O}$, so that (\ref{local average}) holds $\P$-a.s.
We show that $\tilde{w}(x)$ is a Bernoulli random variable under the condition stated in $(i)$, and a deterministic
constant under the condition given in $(ii)$. Let $j\in \N$. On the one hand, the Dominated Convergence Theorem yields that
\begin{equation}\label{conv moments}
\lim_{m\rightarrow \infty}\E\bigg[\bigg(\varphi_m(\R^d)^{-1}\int_{\R^d}\ell(dz)\ \varphi_m(x+z)\bigg)^j\bigg] = \E\big[\tilde{w}(x)^j\big].
\end{equation}

On the other hand, by Fubini's theorem and (\ref{eq duality}), we have that for every $m\in \N$ \setlength\arraycolsep{1pt}
\begin{eqnarray}
\E&\bigg[&\bigg(\varphi_m(\R^d)^{-1}\int_{\R^d}\ell(dz)\ \varphi_m(x+z)\bigg)^j\bigg]\nonumber\\
& & = \varphi_m(\R^d)^{-j}\int_{(\R^d)^j}
dz_1\ldots dz_j \,
\varphi_m(x+z_1)\cdots \varphi_m(x+z_j)\ \E\big[w(t,z_1)\ldots w(t,z_j)\big]
\nonumber\\
& & = \varphi_m(\R^d)^{-j}\int_{(\R^d)^j}
dz_1\ldots dz_j\,
\E\big[w(0,\xi_t^1)\cdots w\big(0,\xi_t^{N_t}\big)\, \big|\,  N_0=j,\ \xi_0^1=z_1,\ldots, \xi_0^j=z_j\big]\times \nonumber\\
& & \qquad \qquad\qquad\qquad\qquad \qquad \qquad \qquad \qquad   \varphi_m(x+z_1)\cdots \varphi_m(x+z_j). \label{dual phi}
\end{eqnarray}
Since the Lebesgue measure of the set of $j$-tuples with at least two identical coordinates is $0$, under the condition of
$(ii)$ the quantity in the right-hand side of (\ref{dual phi}) is equal to
$$
\int_{(\R^d)^j}
dz_1\ldots dz_j \,
 \prod_{i=1}^j \bigg\{\frac{\varphi_m(x+z_i)}{\varphi_m(\R^d)}\,\E_{z_i}[w(0,\xi_t)]\bigg\}
= \bigg(\int_{\R^d}dz\,\frac{\varphi_m(x+z)}{\varphi_m(\R^d)}\,\E_z[w(0,\xi_t)]\bigg)^j.
$$
By our continuity assumption, this quantity tends to $\E_x[w(0,\xi_t)]^j$ as
$m\rightarrow \infty$. Combined with (\ref{conv moments}), this gives us that $\tilde{w}(x)$ is a.s. equal to the
constant $\E_x[w(0,\xi_t)]$ under the condition stated in $(ii)$.

To see $(i)$, consider the case $j=2$ (i.e., $\cA$ consists of two ancestral lineages) and let us write $\tau$ for the
time at which they coalesce, with the convention that $\tau=\infty$ if $\cA$ always contains two lineages.
Since $\varphi_m(x+ \cdot)$ is concentrated on $B(x,2/m)$, using (\ref{cond Bernoulli})
we obtain that for every $\e>0$,
$$
\lim_{m\rightarrow \infty}\frac{1}{\varphi_m(\R^d)^2}\int_{(\R^d)^2}dz_1dz_2\, \varphi_m(x+z_1)\varphi_m(x+z_2)
\P\left[\tau> \e\,|\, N_0=2,\ \xi_0^1=z_1, \xi_0^2=z_2\right] = 0.
$$
Hence, for $j=2$ and $\e<t$, the quantity on the right-hand side of (\ref{dual phi}) can be written
$$
\int_{(\R^d)^2}dz_1dz_2\, \E\big[w(0,\xi_t^1)\ind_{\{\tau\leq \e\}}\, \big|\,  N_0=2,\ \xi_0^1=z_1,\xi_0^2=z_2\big] \frac{\varphi_m(x+z_1)\varphi_m(x+z_2)}{\varphi_m(\R^d)^2}
 + \delta(\e,m),
$$
where $\delta(\e,m)\rightarrow 0$ as $m\rightarrow \infty$ for every fixed $\e$. By the same argument, we have
\begin{eqnarray*}
& & \int_{(\R^d)^2}dz_1dz_2\, \E\big[w(0,\xi_t^1)\ind_{\{\tau\leq \e\}}\,\big|\,  N_0=2,\ \xi_0^1=z_1, \xi_0^2=z_2\big]\frac{\varphi_m(x+z_1)\varphi_m(x+z_2)}{\varphi_m(\R^d)^2}
\\
& &= \int_{(\R^d)^2}dz_1dz_2\, \E\big[w(0,\xi_t^1)\, \big|\,  N_0=2,\ \xi_0^1=z_1,\xi_0^2=z_2\big] \frac{\varphi_m(x+z_1)\varphi_m(x+z_2)}{\varphi_m(\R^d)^2}
+ \delta'(\e,m)\\
& & = \int_{(\R^d)^2}dz_1dz_2\, \E_{z_1}[w(0,\xi_t)] \frac{\varphi_m(x+z_1)\varphi_m(x+z_2)}{\varphi_m(\R^d)^2}  + \delta'(\e,m)\\
& & =  \int_{\R^d}dz_1\, \E_{z_1}[w(0,\xi_t)] \frac{\varphi_m(x+z_1)}{\varphi_m(\R^d)} + \delta'(\e,m),
\end{eqnarray*}
where $\delta'(\e,m)$ also tends to $0$ as $m\rightarrow \infty$ for every $\e>0$, and the third line is justified by
the consistency of $(\cA_t)_{t\geq 0}$.  Using again our continuity assumption on $z\mapsto \E_z[w(0,\xi_t)]$, we obtain
that under the condition stated in $(i)$, the quantity on the right-hand side of (\ref{dual phi}) converges
to $\E_x[w(0,\xi_t)]$ as $m\rightarrow\infty$.  Hence, coming back to (\ref{conv moments}), we arrive at
$$
\E\big[\tilde{w}(x)^2] = \E_x[w(0,\xi_t)]= \E\big[\tilde{w}(x)].
$$
Since $\tilde{w}(x) \in [0,1]$ almost surely, we deduce that $\tilde{w}(x) \in \{0,1\}$ almost surely,
whence $\tilde{w}(x)$ is a Bernoulli random variable.
This completes the proof of Lemma \ref{lem: density} $(i)$. \hfill $\Box$

\medskip

Note that $(ii)$ corroborates a remark at the beginning of \S5 in \cite{EVA1997}. In Evans' construction, all the
genealogical processes used as duals are made up of independent Hunt processes that coalesce instantaneously upon meeting.
Evans points out that, in this case, if $\xi$ and $\xi'$ are two independent processes having the same law as the motion
of a single lineage, then the corresponding $\Xi$-valued process evolves deterministically iff
$$
\mathrm{Vol}\Big(\Big\{(z_1,z_2)\in (\R^d)^2\,:\,\P_{z_1,z_2}\big[\exists\ t\geq 0\,:\, \xi_t=\xi'_t\big]>0\Big\}\Big) = 0.
$$
That is, if the set of pairs of starting points $(z_1,z_2)$ such that $\xi$ and $\xi'$ have a positive chance to meet
in finite time is negligible with respect to Lebesgue measure, then for every $t>0$, $\rho_t$ is a deterministic
function of its initial value (and so is $w(t,\cdot)$). Our proof of Lemma~\ref{lem: density} gives an alternative
proof of Evans's remark when the type-space $K$ is $\{0,1\}$.

\section{Proof of Theorem \ref{theo: A}}\label{section: proof A}

Let us start by proving the convergence stated in Theorem~\ref{theo: A} for a single time $t\geq 0$.
Since we start from $w(0,\cdot)=\ind_H(\cdot)$ (where $H\subset \R^d$ is the half-space of all points whose first
coordinate is non-positive), for every $n\in \N$ we have $w(0,\cdot \sqrt{n})=\ind_H(\cdot)$. Hence, we need only prove the
result for $t>0$.

From our definition of convergence (see (\ref{sense conv})), our aim is to show that for every $j\in \N$ and
$\psi\in C((\R^d)^j)\cap L^1(dx^{\otimes j})$, \setlength\arraycolsep{1pt}
\begin{align*}
\lim_{n\rightarrow \infty}&\E\bigg[\int_{(\R^d)^j}
dx_1\ldots dx_j\,
\psi(x_1,\ldots,x_j) w(tn,x_1\sqrt{n})\cdots w(tn,x_j\sqrt{n})
\bigg] \\
& \qquad \qquad = \E\bigg[\int_{(\R^d)^j}
dx_1\ldots dx_j\,
\psi(x_1,\ldots,x_j)w^{(2)}(t,x_1)\cdots w^{(2)}(t,x_j)
\bigg].
\end{align*}

As we explained in \S\ref{subs: duality}, this question boils down to establishing the asymptotic behaviour of
$$
\int_{(\R^d)^j}
dx_1\ldots dx_j\,
\psi(x_1,\ldots,x_j)\E\big[w(0,\xi_{tn}^1)\cdots w(0,\xi_{tn}^{N_{tn}})\,\big|\, N_0=j,\xi_0^1= x_1\sqrt{n},\ldots,\xi_0^j=x_j\sqrt{n}\big].
$$
This will be achieved in Lemmas~\ref{lem: geneal A1} and~\ref{lem: geneal A2} below,
but first we need some notation.
Recall
that $\xi$ represents the motion of a single ancestral lineage, that is $\xi$ is a compound Poisson process in
which jumps from $x$ to $x+z$ have intensity
$$
m(dz)= \frac{u L_r(z)}{V_r}\ dz.
$$
Observe in passing that this intensity is $0$ whenever $|z|\geq 2r$ (since the start and end points of a jump must
belong to the same ball of radius $r$ and so the size of this jump is bounded by $2r$). For every $n\in \N$, let $\xi^n$
be the process on $\R^d$ defined by
$$
\xi^n_t := \frac{1}{\sqrt{n}}\ \xi_{tn},\qquad t\geq 0,
$$
and let $\cA^n$ be the corresponding rescaling of $\cA$ in which time is multiplied by $n$ and spatial locations
are scaled down by $\sqrt{n}$. More formally, here we view $\cA^n$ and $\cA$ as having values in the state-space $\bigcup_{m \ge 0} \{m \} \times (\R^d)^m$; thus the first coordinate of $\cA$ or $\cA_n$ indicates the number of distinct lineages, while the remaining coordinates give their respective positions. (In Section \ref{S:heavy} we will find it convenient to enrich the state space to also record the whole genealogical information).

\begin{lemma}\label{lem: geneal A1}
If $d=1$, for every $j\in \N$ and $x_1,\ldots,x_j\in \R^d$, the process $\cA^n$ starting from $j$ lineages at
locations $x_1,\ldots,x_j$ converges, in the sense of finite-dimensional distributions, to a system $\cA^{\infty}$ of
independent Brownian motions with clock speed $\sigma^2$ that coalesce instantaneously upon meeting.

More generally, let $k\in\N$ and $0<t_1<\ldots<t_k$.  Suppose that we start $\cA^n$ with $j_0$ lineages at distinct
locations $x_{0,1},\ldots,x_{0,j_0}$, let the process evolve until time $t_1$, add to the surviving lineages $j_1$
lineages at distinct locations $x_{1,1},\ldots,x_{1,j_1}$, let all resulting lineages evolve until time $t_2$ when we add $j_2$
further lineages, and so on.  Call the corresponding process $\hat{\cA}^n$.  Define $\hat{\cA}^{\infty}$ analogously.
Then for any $t\geq 0$, the law of $\hat{\cA}^n_t$ converges to that of $\hat{\cA}^{\infty}_t$ as $n$ tends to infinity.
\end{lemma}

\begin{lemma}\label{lem: geneal A2}
If $d\geq 2$, for every $j\in \N$ and distinct $x_1,\ldots,x_j\in \R^d$, the process $\cA^n$ starting from $j$ lineages at
locations $x_1,\ldots,x_j$ converges to a system of independent Brownian motions with speed $\sigma^2$.
In particular, the limiting lineages never coalesce.

More generally, define $\hat{\cA}^n$ and $\hat{\cA}^{\infty}$ as in Lemma~\ref{lem: geneal A1}.
Then for any $t\geq 0$, the law of $\hat{\cA}^n_t$ converges to that of $\hat{\cA}^{\infty}_t$ as $n$ tends to infinity.
\end{lemma}
We postpone the proofs of Lemmas~\ref{lem: geneal A1} and \ref{lem: geneal A2} until the end of this section.

Since the boundary of $H$ has zero Lebesgue measure, Portmanteau's Lemma and the first part of
Lemma~\ref{lem: geneal A1} give us that if $d=1$, (using the obvious generalisation to $\cA^{\infty}$ of our previous notation)
\begin{eqnarray*}
& & \lim_{n\rightarrow \infty}\int_{(\R^d)^j}
dx_1\ldots dx_j \,
\psi(x_1,\ldots,x_j)\E\big[\ind_H(\xi_{tn}^1)\cdots \ind_H\big(\xi_{tn}^{N_{tn}}\big)\,\big|\, N_0=j,\xi_0^1= x_1\sqrt{n},\ldots,\xi_0^j=x_j\sqrt{n}\big]
\\
& & = \lim_{n\rightarrow \infty}\int_{(\R^d)^j}
dx_1\ldots dx_j \,
\psi(x_1,\ldots,x_j)\E\big[\ind_H(\xi_t^{n,1})\cdots \ind_H\big(\xi_t^{n,N^n_t}\big)\big| N^n_0=j,\xi_0^{n,1}= x_1,\ldots,\xi_0^{n,j}=x_j\big]
\\
& & = \int_{(\R^d)^j}
dx_1\ldots dx_j\,
\psi(x_1,\ldots,x_j)\E\big[\ind_H(\xi_t^{\infty,1})\cdots \ind_H\big(\xi_t^{\infty,N^{\infty}_t}\big)\,\big|\, N_0^{\infty}=j,\xi_0^{\infty,1}= x_1,\ldots,\xi_0^{\infty,j}=x_j\big].
\end{eqnarray*}
Now, Theorem~$4.1$ in \cite{EVA1997} guarantees that there exists a unique $\Xi$-valued Markov process
starting from (the equivalence class of) $\ind_H(x)$ and dual to $\cA^{\infty}$ through the relations (\ref{eq duality}).
Let us call this process $\rho^{(2)}$.  Using the more compact notation of \S\ref{subs: form conv}, we obtain that for
every $j\in\N$ and $\psi\in C((\R^d)^j)\cap L^1(dx^{\otimes j})$,
$$
\lim_{n\rightarrow \infty}\E\big[I_j(\rho_t^n;\,\psi)\big]= \E\big[I_j(\rho_t^{(2)};\,\psi)\big].
$$
Since this family of test functions in dense in $C(\Xi)$ (c.f.~\S\ref{subs: form conv}),
we can conclude that $\rho^n_t\stackrel{\mathcal{L}}{\rightarrow}\rho^{(2)}_t$ as $n\rightarrow \infty$.
It is then straightforward to check that the conditions of Lemma~\ref{lem: density} $(i)$ are satisfied,
and so for a.e. $x\in \R^d$, $w^{(2)}(t,x)$ is a Bernoulli random variable with parameter
$$
\P_x\big[\xi_t^{\infty}\in H\big] = \P_x\big[X_{\sigma^2t}\in H\big] = p^2(t,x).
$$
Moreover, by Lemma~\ref{lem: geneal A1} and (\ref{point duality}), the correlations between the values of $w^{(2)}(t,\cdot)$
at different sites can be described as follows.  For every $j\in \N$ and Lebesgue-a.e. $(x_1,\ldots,x_j)$, \setlength\arraycolsep{1pt}
\begin{eqnarray}
\E\big[&w^{(2)}&(t,x_1)\ldots w^{(2)}(t,x_j)\,\big] \nonumber \\
  & & = \E\Big[w^{(2)}(0,\xi_t^{\infty,1})\cdots w^{(2)}\big(0,\xi_t^{\infty,N^{\infty}_t}\big)\, \Big|\, N^{\infty}_0=j,\ \xi_0^{\infty,1}=x_1,\ldots, \xi_0^{\infty,j}=x_j\Big]\nonumber \\
& & =  \P\Big[\xi_t^{\infty,i}\in H,\ \forall i\in \{1,\ldots,N^{\infty}_t\} \Big|\, N^{\infty}_0=j,\ \xi_0^{\infty,1}=x_1,\ldots, \xi_0^{\infty,j}=x_j\Big]. \label{correl A1}
\end{eqnarray}
Since we are dealing with Bernoulli random variables, equation~(\ref{correl A1}) completely characterizes these correlations.

If $d\geq 2$, by the same chain of arguments (using this time Lemma~\ref{lem: geneal A2}), we obtain
\begin{eqnarray*}
& & \lim_{n\rightarrow \infty}\int_{(\R^d)^j}
dx_1\ldots dx_j \,
\psi(x_1,\ldots,x_j)\E\big[\ind_H(\xi_{tn}^1)\cdots \ind_H(\xi_{tn}^{N_{tn}})\,\big|\, N_0=j,\xi_0^1= x_1\sqrt{n},\ldots,\xi_0^j=x_j\sqrt{n}\big]
\\
& &\qquad = \int_{(\R^d)^j}
dx_1\ldots dx_j \,
\psi(x_1,\ldots,x_j)\P_{x_1}\big[\xi_t^{\infty,1}\in H\big]\cdots \P_{x_j}\big[\xi_t^{\infty,j}\in H\big].
\end{eqnarray*}
Here again, these equalities guarantee the convergence in law of $\rho^n_t$ towards the value at time $t$ of the
unique $\Xi$-valued Markov process $\rho^{(2)}$ starting from $\ind_H(x)$ and dual to the system $\cA^{\infty}$ of
independent Brownian motions which never coalesce. Lemma~\ref{lem: density} $(ii)$ then applies and gives us that
for a.e. $x\in \R^d$, $w^{(2)}(t,x)$ is the deterministic constant $p^2(t,x)$.

So far, we have obtained the desired convergence at a given time $t>0$, and the form of the local densities of $1$'s in the limit.
It remains to show that the convergence holds true for finitely many times $0\leq t_1\leq \cdots \leq t_k$. Because functions of the form $I_j(\cdot\,;\,\psi)$ are dense in $C(\Xi)$, we need only show that for every $j_1,\ldots,j_k$ and $\psi_1,\ldots, \psi_k$,
\begin{equation}\label{conv fd}
\lim_{n\rightarrow \infty}\E\bigg[\prod_{i=1}^k I_{j_i}\big(\rho^n_{t_i}\, ;\,\psi_i\big)\bigg] = \E\bigg[\prod_{i=1}^k I_{j_i}\big(\rho^{(2)}_{t_i}\,;\, \psi_i\big)\bigg].
\end{equation}
Therefore, let us fix $j_1,\ldots,j_k$ and $\psi_1, \ldots,\psi_k$ such that $\psi_i \in C((\R^d)^{j_i})\cap L^1(dx^{\otimes j_i})$.
To simplify notation, we write $\und{x}^i$ for the vector $(x^i_1,\ldots,x^i_{j_i})$ and $W_i^n(\und{x}^i)$ for the
product $\prod_{l=1}^{j_i}w^n(t_i,x^i_l)$.  We will occasionally abuse notations and write $\und x^i \cup \und x^j$ for the concatenation of the vectors $\und x^i $ and $\und x^j$.
Our strategy is to use duality again, but now with the genealogical process
described in the second part of Lemmas~\ref{lem: geneal A1} and \ref{lem: geneal A2}.  Once again, to simplify our notation,
let us denote the law of $\cA^n$ (resp., $\cA^n_t$) starting from $j$ lineages at locations $\und{x}=(x_1,\ldots,x_j)$
by $\P^n_{\und{x}}$ (resp., $\P^n_{\und{x},t}$). Hence recall that $\P^n_{\und x}$ is a distribution over $\bigcup_{m \ge 0} \{m \} \times (\R^d)^m$. Using the Markov property of $w$ at time $t_{k-1}n$ and the duality property~(\ref{eq duality}), we can write
\setlength\arraycolsep{1pt}
\begin{eqnarray*}
\E&\bigg[&\prod_{i=1}^k I_{j_i}\big(\rho^n_{t_i}\,;\,\psi_i\big)\bigg] \\
& & = \int\ldots\int d\und{x}^1\cdots d\und{x}^k \ \psi_1(\und{x}^1)\cdots \psi_k(\und{x}^k) \\
& & \qquad \qquad \qquad \qquad\times \E\bigg[\bigg\{\prod_{i=1}^{k-1}W^n_i\big(\und{x}^i\big)\bigg\}\ \E^n_{\und{x}^k}\Big[w^n\big(t_{k-1},\xi^{n,1}_{t_k-t_{k-1}}\big)\cdots\, w^n\Big(t_{k-1},\xi^{n,N^n_{t_k-t_{k-1}}}_{t_k-t_{k-1}}\Big)\Big]\bigg]\\
& & = \int\ldots\int d\und{x}^1\cdots d\und{x}^k \ \psi_1(\und{x}^1)\cdots \psi_k(\und{x}^k) \int d\P^n_{\und{x}^k,t_k-t_{k-1}}\big(m_{k-1},y^{k-1}_1,\ldots,y^{k-1}_{m_{k-1}}\big)\\
& & \quad  \E\bigg[\bigg\{\prod_{i=1}^{k-2}W^n_i\big(\und{x}^i\big)\bigg\}\ w^n\big(t_{k-1},x^{k-1}_1\big)\cdots w^n\big(t_{k-1},x^{k-1}_{j_{k-1}}\big) w^n\big(t_{k-1},y^{k-1}_1\big)\cdots\ w^n\big(t_{k-1},y^{k-1}_{m_{k-1}}\big)\bigg].
\end{eqnarray*}
Since the law of the locations at time $t_k-t_{k-1}$ of the $N^n_{t_k-t_{k-1}}$ lineages is absolutely continuous with
respect to Lebesgue measure, we can carry on the recursion and use the Markov property (this time at time $t_{k-2}$) and duality to
write the quantity above as \setlength\arraycolsep{1pt}
\begin{eqnarray}
\int\ldots\int& & d\und{x}^1\cdots d\und{x}^k \ \psi_1(\und{x}^1)\cdots \psi_k(\und{x}^k) \int d\P^n_{\und{x}^k,t_k-t_{k-1}}\big(m_{k-1},y^{k-1}_1,\ldots,y^{k-1}_{m_{k-1}}\big) \nonumber\\
\times \int& & d\P^n_{\und{x}^{k-1}\cup \und{y}^{k-1},t_{k-1}-t_{k-2}}\big(m_{k-2},y^{k-2}_1,\ldots,y^{k-2}_{m_{k-2}}\big)\E\bigg[\bigg\{\prod_{i=1}^{k-3}W^n_i\big(\und{x}^i\big)\bigg\} \nonumber\\
& & \qquad \qquad \times w^n\big(t_{k-2},x^{k-2}_1\big)\cdots\ w^n\big(t_{k-2},x^{k-2}_{j_{k-2}}\big) w^n\big(t_{k-2},y^{k-2}_1\big)\cdots\ w^n\big(t_{k-2},y^{k-2}_{m_{k-2}}\big)\bigg]\nonumber\\
= \int\ldots\int & &  d\und{x}^1\cdots d\und{x}^k \ \psi_1(\und{x}^1)\cdots \psi_k(\und{x}^k) \int d\P^n_{\und{x}^k,t_k-t_{k-1}}\big(m_{k-1},y^{k-1}_1,\ldots,y^{k-1}_{m_{k-1}}\big) \nonumber\\
\times \int & &\cdots \int d\P^n_{\und{x}^1\cup \und{y}^1,t_1}\big(m_0,y^0_1,\ldots,y^0_{m_0}\big)\E\Big[w^n\big(0,y^0_1) \cdots w^n\big(0,y^0_{m_0}\big)\Big].\label{multiple dual}
\end{eqnarray}
Now, recall the family of processes $\hat{\cA}^n$ introduced in the second part of Lemmas~\ref{lem: geneal A1} and \ref{lem: geneal A2}.
Let us denote the times of appearance and the locations of the additional lineages in the form
$(\tau_1,\und{z}^1),\ldots, (\tau_k,\und{z}^k)$.  Using (recursively) the Markov property of $\hat{\cA}^n$, we obtain that the
quantity on the right-hand side of (\ref{multiple dual}) is equal to
\begin{eqnarray*}
\int&\ldots&\int d\und{x}^1\cdots d\und{x}^k \ \psi_1(\und{x}^1)\cdots \psi_k(\und{x}^k)\\
& & \qquad \qquad \times\E\Big[w^n\big(0,\hat{\xi}_{t_k}^{\ n,1}\big)\cdots \ w^n\Big(0,\hat{\xi}_{t_k}^{\ n,\hat{N}^n_{t_k}}\Big)\,\Big|\,\big(0,\und{x}^k\big),\big(t_k-t_{k-1},\und{x}^{k-1}\big),\ldots,\big(t_k-t_1,\und{x}^1\big)\Big].
\end{eqnarray*}
Let us now conclude when $d=1$ (the reasoning is exactly the same when $d\geq 2$).  Recall that for every
$n\in \N$, $w^n(0,\cdot)=\ind_H(\cdot)= w^{(2)}(0,\cdot)$. By the second part of Lemma~\ref{lem: geneal A1} and
the Dominated Convergence Theorem (and the fact that the boundary of $H$ has zero Lebesgue measure), we obtain that
\begin{eqnarray*}
\lim_{n\rightarrow \infty}& \E&\bigg[\prod_{i=1}^k I_{j_i}\big(\rho^n_{t_i}\, ;\,\psi_i\big)\bigg] \\
& = & \int\ldots\int d\und{x}^1\cdots d\und{x}^k \ \psi_1(\und{x}^1)\cdots \psi_k(\und{x}^k)\\
& & \qquad \ \times\E\Big[w^{(2)}\big(0,\hat{\xi}_{t_k}^{\ \infty,1}\big)\cdots w^{(2)}\Big(0,\hat{\xi}_{t_k}^{\ \infty,\hat{N}^{\infty}_{t_k}}\Big)\,\Big|\,\big(0,\und{x}^k\big),\big(t_k-t_{k-1},\und{x}^{k-1}\big),\ldots,\big(t_k-t_1,\und{x}^1\big)\Big].
\end{eqnarray*}
Analogous calculations using the duality between $\cA^{\infty}$ and $\rho^{(2)}$ lead to (\ref{conv fd}).
This completes the proof of Theorem~\ref{theo: A}. \hfill $\Box$

\medskip
It remains to prove Lemmas~\ref{lem: geneal A1} and \ref{lem: geneal A2}. Let us start with the latter, which is
somewhat simpler, but contains the main ingredients of both proofs.

\medskip
\noindent{\bf Proof of Lemma \ref{lem: geneal A2}.} Let $x_1,\ldots,x_k$ be $k$ distinct points of $\R^d$.
Suppose that $\cA^n$ starts from $k$ lineages at locations $x_1\sqrt{n},\ldots,x_k\sqrt{n}$. First, since a single lineage $\xi$
follows a finite-rate homogeneous jump process whose jumps are uniformly bounded by $2r$, standard arguments guarantee
that $\xi^n = (n^{-1/2}\xi_{tn})_{t\geq 0}$ converges in distribution to Brownian motion with clock speed $\sigma^2$ given in (\ref{def sigma}).

Second, observe that two lineages can be hit by the same event (and possibly coalesce) only if they lie at distance at most $2r$ of
each other. Consequently, as long as they are at distance greater than $2r$ they evolve independently, according to the law of the
motion of a single lineage. Hence, let us define $n\tau_n$ to be the first time at which at least two of the $k$ initial lineages
are within distance at most $2r$ of one another.  Equivalently, $\tau_n$ is the first time at which at least two lineages of
$\cA^n$ are at separation at most $2r/\sqrt{n}$. We wish to show that for any $t\geq 0$, $\P^n_{\und{x}}[\tau_n\leq t]\rightarrow 0$
as $n\rightarrow \infty$.

To this end, note that until time $\tau_n$, the motions of the rescaled lineages $\xi^{n,1}, \ldots, \xi^{n,k}$ can be
embedded in the paths of independent standard Brownian motions $X^1,\ldots,X^k$ starting from $x_1,\ldots,x_k$ (we use the
same Brownian motions for all $n$). Indeed, for each path $i$ we proceed as follows (this construction is in the spirit of
the one-dimensional Skorokhod Embedding Theorem, see e.g. \cite{BIL1995}).  Let $(R^{n,i}_j)_{j\geq 1}$ be a sequence of i.i.d.
random variables (independent of $X^i$) distributed according to the law of the radius of a typical jump of $\xi^n$, and let
us define a sequence $\{s_{i,j}^n,j\geq 0\}$ of random times, recursively, by
\begin{enumerate}
\item $s^n_{i,0}:=0$,
\item for every $j\geq 1$, $s^n_{i,j}$ is the first time greater than $s^n_{i,j-1}$ at which $X^i$ exits the ball $B\big(X^i_{s^n_{i,j-1}}, R^{n,i}_j\big)$.
\end{enumerate}
By rotational symmetry of the law of a jump of $\xi^{n,i}$, conditional on its radius being $\gamma$ the location of $\xi^{n,i}$
just after the jump is uniformly distributed over the sphere $\partial B(\xi^{n,i}_{t-},\gamma)$.
Likewise, conditional on the variable $R^{n,i}_j$ being equal to $\gamma$, the location of $X^i_{s^n_{i,j}}$ is uniformly
distributed over $\partial B(X^i_{s^n_{i,j-1}}, \gamma)$.  Consequently, by comparing their jump rates and their jump distributions,
one can show that for every $i\in \{1,\ldots,k\}$ the processes $(\xi^{n,i}_t)_{t\geq 0}$ and
$\big(X^i_{s_{i,j(n,i,t)}^n}\big)_{t\geq 0}$ have the same laws, where $(j(n,i,t))_{t\geq 0}$ is a Poisson process
with intensity $n u V_r$ (recall from (\ref{cond existence}) that $uV_r$ is the jump rate of an unrescaled lineage
under the conditions of Case A, where $V_r$ is the volume of a ball of radius $r$).
Since the lineages $\xi^{n,1},\ldots,\xi^{n,j}$ evolve independently until time $\tau_n$, we can ask that the
Poisson processes $\{j(n,1,\cdot),\ldots,j(n,k,\cdot)\}$ should be independent and the embedding holds for
all $i\in \{1,\ldots,k\}$ simultaneously until the first time $t$ such that
$$
\big|X^i_{s_{i,j(n,i,t)}^n}- X^m_{s_{m,j(n,m,t)}^n}\big|\leq 2r/\sqrt{n}\qquad \mathrm{for\ some}\ i\neq m.
$$

Now, each rescaled lineage makes jumps of size at most $2rn^{-1/2}$ at rate $\cO(n)$.  Hence, each difference
$s^n_{i,j}-s^n_{i,j-1}$ is the exit time of Brownian motion from a ball of radius $\cO(n^{-1/2})$, and
$s_{i,j(n,i,t\wedge \tau_n)}^n$ is the sum of (morally) $\cO(n)$ such quantities, all independent of one another.
More formally, if we write $R$ for the (random) radius of a typical jump of an unrescaled lineage and if we
notice that the exit time of Brownian motion starting at $0$ from a ball $B(0,\gamma)$ is bounded by the first time
that one of its coordinates leaves the interval $[-\gamma,\gamma]$, then for all $n\in \N$ and all $1\leq i\leq k$
we can write
$$
\E\big[s^{n}_{i,1}\big]\leq d\, \E\Big[\big(R^{n,i}_1\big)^2\Big] = \frac{d}{n}\,\E\big[R^2\big]\leq \frac{4dr^2}{n},
$$
where the first inequality uses the property that the exit time from $[-\gamma,\gamma]$ of one-dimensional
Brownian motion starting at $0$ has expectation $\gamma^2$. By the independence of $X_i$ and the Poisson processes,
this yields that for all $n$ and $i$,
$$
\E\big[s_{i,j(n,i,t\wedge \tau_n)}^n\big] = \E\big[j(n,i,t\wedge \tau_n)\big].\E\big[s^n_{i,1}\big]\leq 4duV_rr^2.
$$

To conclude our proof, let us observe that $\P^n_{\und{x}}[\tau_n \leq t]$ is bounded by the probability that at least
two of the $k$ independent Brownian motions $X^1,\ldots,X^k$ come within distance $2rn^{-1/2}$ before
time $\min\{s^n_{i,j(n,i,t)},\, 1\leq i\leq k\}$. But if $\tilde{\tau}_n$ denotes the first time at which two
independent Brownian motions starting at $x_1\neq x_2$ come within distance $2rn^{-1/2}$, for every $T\geq 0$ we have
$$
\lim_{n\rightarrow \infty} \P_{x_1,x_2}\big[\tilde{\tau}_n \leq T\big] = 0.
$$
Hence, the probability that at least two out of $k$ independent Brownian motions come within distance $2rn^{-1/2}$
before any given time $T$ also tends to $0$, and thanks to the uniform bound on the expectation of
$s_{i,j(n,i,t\wedge \tau_n)}^n$ (together with the Markov inequality), it is straightforward to obtain that for any $t\geq 0$
$$
\lim_{n\rightarrow \infty} \P^n_{\und{x}}\big[\tau_n \leq t\big] = 0.
$$
We have thus shown that with probability growing to $1$ as $n\rightarrow \infty$, until a given time $t\geq 0$ the $k$
ancestral lineages evolve as if they were independent. Since the law of each $\xi^{n,i}$ converges to that of
Brownian motion with clock speed $\sigma^2$, the convergence of the one-dimensional distributions of $\cA^n$ to those of a
collection of $k$ independent Brownian motions is proved.

The proofs of the convergence of the finite-dimensional distributions and that of the second part of Lemma~\ref{lem: geneal A2}
follow the same lines, using the Markov property of each $\cA^n$ at suitable times.  Details are left to the reader.\hfill $\Box$

\medskip
\noindent{\bf Proof of Lemma \ref{lem: geneal A1}.}
Once again we start with the one-dimensional distributions, and proceed by recursion on the number $m$ of lineages of $\cA^n$.
As in the proof of Lemma~\ref{lem: geneal A2}, before rescaling each lineage follows a homogeneous symmetric (finite rate)
jump process, whose jumps have length at most $2r$, and so $\xi^n=(n^{-1/2}\xi_{nt})_{t\geq 0}$ converges in distribution to
Brownian motion with clock speed $\sigma^2$ as $n$ tends to infinity.

Let us consider the case $m=2$.  As we saw in the proof of Lemma~\ref{lem: geneal A2}, the two rescaled lineages evolve
independently until they come within distance $2rn^{-1/2}$ of one another.  Let us first show that this `meeting' time converges to
the meeting time (at distance $0$) of two independent Brownian motions starting at $x_1$ and $x_2$ and with clock speed $\sigma^2$,
and secondly that coalescence is quasi-instantaneous once the lineages are gathered at this distance.

For the first claim, let us write $\tau_n$ for the time at which $\xi^{n,1}$ and $\xi^{n,2}$ first come within
distance at most $4rn^{-1/2}$ of one another
(note the constant $4$ instead of $2$, which we shall need later for purely technical reasons).
Because the motion of a single lineage is a symmetric jump process, until $\tau_n$ the law of the
difference $\xi^{n,1}-\xi^{n,2}$ is the same as that of the motion of a single rescaled lineage, run at speed $2$.
Let $X$ be a standard one-dimensional Brownian motion, starting from $x_1-x_2$ and independent of all $\xi^n$'s.
Using anew the construction introduced in the proof of Lemma~\ref{lem: geneal A2}, for every $n$ we can find a sequence of
random times $\{s^n_j, j\geq 0\}$ such that $(\xi^{n,1}_t-\xi^{n,2}_t)_{t\geq 0}$ has the same law as $(X_{s^n_{j(n,t)}})_{t\geq 0}$,
where $j(n,\cdot)$ is a Poisson process, independent of $X$ and with intensity $4nru$ (that is, twice the jump rate of a
single rescaled lineage). Recall from the proof of Lemma~\ref{lem: geneal A2} that for every $n\in \N$, the random
variables $s^n_j-s^n_{j-1}$, $j\geq 1$, are i.i.d and if $R$ is distributed like the radius of a typical jump of $\xi$,
we have $\E[ns^n_1]=\E[R^2]<\infty$. Let $t\geq 0$, and, as a first step, let us show that
$s^n_{j(n,t)}$ converges in probability towards $2\sigma^2t$ as $n$ grows to infinity. The second step will then consist
of proving that, for every $t\geq 0$, the probability that $\tau_n >t$ tends to the probability that the hitting time of
$0$ by $X$ is greater than $2\sigma^2t$. This will give us the desired result.

By definition, $j(n,t)$ is a Poisson random variable with parameter $(4nur)t$. By the Central Limit Theorem, we
therefore have that
\begin{equation}\label{conv meeting 1}
n^{-1/2}\big(j(n,t)- 4nurt\big) \stackrel{(d)}{\longrightarrow} \mathcal{N}(0,4urt).
\end{equation}
Now, recalling the properties of the $s^n_i-s^n_{i-1}$'s expounded above, by the Strong Law of Large Numbers we have
\begin{equation}\label{conv meeting 2}
s^n_{\lfloor 4nurt \rfloor} = \frac{1}{n}\sum_{i=1}^{\lfloor 4nurt\rfloor} n\big(s^n_i - s^n_{i-1}\big) \stackrel{\mathrm{a.s.}}{\longrightarrow} 4urt \times \E[R^2] \qquad \mathrm{as}\ n\rightarrow \infty,
\end{equation}
where $\lfloor z \rfloor$ denotes the integer part of $z$. But $\sigma^2$ is defined in (\ref{def sigma}) as the variance of the displacement at time $1$ of a single unrescaled lineage, and so
$$
\sigma^2 = 2ur\,\E[R^2],
$$
which shows that the limit in (\ref{conv meeting 2}) is equal to $2\sigma^2t$. To conclude the first step, observe that $|s^n_{j(n,t)}-s^n_{\lfloor 4nurt \rfloor}|$ is the sum of $|j(n,t)- \lfloor 4nurt \rfloor|$ i.i.d. terms of the form $s^n_i-s^n_{i-1}$, all of them independent of $j(n,t)$, so that for every $\e>0$ and every $n\geq 1$ we have
$$
\P\big[\big|s^n_{j(n,t)}-s^n_{\lfloor 4nurt \rfloor}\big|>\e\big] \leq  \P\big[|j(n,t)- 4nurt|> n^{3/4}\big] + \P\left[\sum_{i=1}^{n^{3/4}}\big(s^n_i-s^n_{i-1}\big)>\e\right].
$$
As $n\rightarrow \infty$, the first term on the right-hand side tends to $0$ by (\ref{conv meeting 1}), while Markov's inequality gives us that
$$
\P\left[\sum_{i=1}^{n^{3/4}}\big(s^n_i-s^n_{i-1}\big)>\e\right] \leq \frac{1}{\e}\ \E\left[\sum_{i=1}^{n^{3/4}}\big(s^n_i-s^n_{i-1}\big)\right] = \frac{C}{\e n^{1/4}} \longrightarrow 0.
$$
Since this is true for any $\e>0$, $s^n_{j(n,t)}-s^n_{\lfloor 4nurt\rfloor}$ converges in probability to $0$.
But we have shown that $s^n_{\lfloor 4nurt\rfloor}$ converges a.s. to $2\sigma^2t$, and so we obtain that
$s^n_{j(n,t)}$ converges in probability to $2\sigma^2t$, as required.

As explained above, we can now use this result to show that $\tau_n$ converges in distribution to the hitting time of
$0$ by $(X_{2\sigma^2t})_{t\geq 0}$. Indeed, by construction of the random times $s^n_i$ and the fact that the rescaled jumps
of a lineage are bounded by $2r/\sqrt{n}$, for any $i\geq 1$ the Brownian motion $X$ cannot move to a distance greater
than $2r/\sqrt{n}$ from $X_{s^n_{i-1}}$ before time $s^n_i$. Thus, if $\tau_0$ denotes the hitting time of $0$ by $X$, we have
$$
\P_{x_1-x_2}[\tau_n>t] \leq \P_{x_1-x_2}\big[\tau_0> s^n_{j(n,t)}\big].
$$
But we showed that $s^n_{j(n,t)}$ converges in probability towards $2\sigma^2t$ as $n\rightarrow \infty$, and so
\begin{equation}\label{conv meeting 3}
\limsup_{n\rightarrow \infty}\ \P_{x_1-x_2}[\tau_n>t] \leq \P_{x_1-x_2}\big[\tau_0> 2\sigma^2t\big].
\end{equation}
On the other hand, for every $\e\in (0,|x_1-x_2|/2)$ and every $n$ large enough, we can write
\begin{eqnarray*}
\P_{x_1-x_2}[\tau_n>t] &\geq &\P_{x_1-x_2}\big[X\ \mathrm{does\ not\ enter\ }B(0,4r/\sqrt{n})\ \mathrm{before}\ s^n_{j(n,t)}\big]\\
& \geq & \P_{x_1-x_2}\big[X\ \mathrm{does\ not\ enter\ }B(0,\e)\ \mathrm{before}\ s^n_{j(n,t)}\big].
\end{eqnarray*}
Again, we can deduce from the convergence in probability of $s^n_{j(n,t)}$ to $2\sigma^2t$ that
$$
\liminf_{n\rightarrow \infty}\ \P_{x_1-x_2}[\tau_n>t] \geq \P_{x_1-x_2}\big[X\ \mathrm{does\ not\ enter\ }B(0,\e)\ \mathrm{before}\ 2\sigma^2t\big].
$$
This inequality holds for every small $\e>0$, and by the point recurrence of one-dimensional Brownian motion, we can conclude that
\begin{equation}\label{conv meeting 4}
\liminf_{n\rightarrow \infty}\ \P_{x_1-x_2}[\tau_n>t] \geq \P_{x_1-x_2}\big[\tau_0> 2\sigma^2t\big].
\end{equation}
Together with (\ref{conv meeting 3}), we obtain that for every $t>0$
\begin{equation}\label{conv meeting}
\lim_{n\rightarrow \infty} \P_{x_1-x_2}[\tau_n>t] = \P_{x_1-x_2}\big[\tau_0> 2\sigma^2t\big],
\end{equation}
from which we can conclude that the `meeting time at distance $4r/\sqrt{n}$' of two rescaled lineages
starting at $x_1$ and $x_2$ converges in distribution to the hitting time of $0$ by Brownian motion with clock speed $2\sigma^2$,
or equivalently to the meeting time of 2 independent Brownian motions each of clock speed $\sigma^2$.

Let us now prove our second claim; that is, let us show that once at distance at most $4r/\sqrt{n}$, the additional time
the two lineages need to merge becomes negligible as $n$ tends to infinity. Because the proof is highly reminiscent of
that of Proposition 6.4$(b)$ in \cite{BEV2010}, we only outline the main steps here.  Let us work with the unrescaled lineages,
and suppose they start at distance at most $4r$ of each other. First, it is not difficult to convince oneself that the
first time at which the two lineages are at separation less than $2r$ is of order $\cO(1)$, `uniformly' over all initial
locations which are at separation at most $4r$. Once close together, they become correlated, because they can be hit by the
same reproduction event. But for the same reason, they have a positive probability of being affected by the same event and of
coalescing before separating again to distance at least $2r$. If they do coalesce, the additional time they had to wait for this
event is also of order $\cO(1)$. If they separate rather than coalescing, then again the time they need to come back
to separation less than $2r$ is of order $\cO(1)$, and once `gathered' they have a positive chance to coalesce before separating,
and so on. In the end, the number of excursions of $\xi^1-\xi^2$ out of $B(0,2r)$ before the two lineages merge can be
stochastically bounded by a geometric random variable, and each of the finitely many excursions and incursions lasts a time
of order $\cO(1)$. This tells us that for every $\e>0$, one can find $T(\e)>0$ such that
$$
\sup_{|y_1-y_2|\leq 4r}\ \P_{(y_1,y_2)}\big[\xi^1\ \mathrm{and}\ \xi^2\ \mathrm{do\ not\ coalesce\ before}\ T(\e)\big] \leq \e.
$$
Rephrasing the above inequality in terms of the rescaled lineages, we obtain that, for every $n\geq 1$,
\begin{equation}\label{instant coal}
\sup_{|z_1-z_2|\leq 4r/\sqrt{n}}\ \P_{(z_1,z_2)}\big[\xi^{n,1}\ \mathrm{and}\ \xi^{n,2}\ \mathrm{do\ not\ coalesce\ before}\ T(\e)/n\big] \leq \e.
\end{equation}
Finally, if $\tau^c_n$ denotes the coalescence time of $\xi^{n,1}$ and $\xi^{n,2}$, using the strong Markov property
of $(\xi^{n,1},\xi^{n,2})$ at time $\tau_n$, we have, for every $t>0$,
$$
\P_{(x_1,x_2)}[\tau^c_n-\tau_n >t]  = \E_{(x_1,x_2)}\Big[\ind_{\{\tau_n<\infty\}}\P_{(\xi^{n,1}_{\tau_n},\xi^{n,2}_{\tau_n})}[\tau^c_n>t]\Big].
$$
By (\ref{instant coal}), the probability inside the expectation tends to $0$ as $n\rightarrow \infty$, and so does the quantity on
the left-hand side (by dominated convergence).  Hence, $\tau_n^c-\tau_n$ converges to $0$ in probability. This concludes the proof of
the first part of Lemma \ref{lem: geneal A1} when $m=2$: in the limit, the two lineages follow independent Brownian motions run at
clock speed $\sigma^2$ until the first time at which they meet, which is also the time at which they coalesce by the convergence
of $\tau_n^c-\tau_n$ to $0$.

We now proceed by induction.
Suppose we know that the result of Lemma \ref{lem: geneal A1} holds true for a system of $m-1$ lineages.
Let $x_1<\ldots<x_m$ be $m$ distinct points of $\R$ and suppose that $m$ lineages start from these locations.
Because the lineages `choose' to take part in an event that encompasses them independently of one another, the law of the
restriction of the system started from $m$ lineages to that started from $m-1$ lineages at $x_1,\ldots,x_{m-1}$ is the same as that
of the $(m-1)$-system starting from $x_1,\ldots,x_{m-1}$. (This is the `consistency' of the genealogical process described below
Lemma~\ref{lem: density}).  Hence, our inductive hypothesis tells us that the restricted process converges to a system of
(initially) $m-1$ independent Brownian motions with clock speed $\sigma^2$, that coalesce instantaneously upon meeting.
Now, as we explained several times already, the motion of the $m$-th lineage, starting at the right-most location $x_m$, is independent
of  that of the others until the first time, $\tau_n$, at which it comes to within distance $2r/\sqrt{n}$ of another lineage.
But with probability tending to $1$, the right-most lineage among those that started from $x_1,\ldots,x_{m-1}$ is the lineage ancestral
to the individual sampled in $x_{m-1}$. Indeed, our inductive hypothesis guarantees that the probability that the lineage starting
from $x_{m-1}$ jumps over a lineage on its left without coalescing with it tends to $0$ as $n$ tends to infinity.
Again by consistency of the genealogical process, when singled out, the motion of lineage $m-1$ has the same law as the process
$\xi^n$ (that is, a typical single lineage), and so we can focus on the two right-most lineages and use the results obtained for
$m=2$ to conclude: their meeting time at distance at most $4r/\sqrt{n}$ converges in distribution to the meeting time of two
independent Brownian motions run at clock speed $\sigma^2$, and in the limit this meeting time is also the coalescence time of the
two lineages. But this is precisely the evolution of a system of (initially) $m$ independent Brownian motions which coalesce
instantaneously when they meet, and so the desired convergence also holds for a system starting with $m$ lineages.

As in the proof of Lemma \ref{lem: geneal A2}, the other points of Lemma \ref{lem: geneal A1} are obtained by using the convergence of the one-dimensional distributions and the Markov property at suitable times. \hfill $\Box$

\section{Heavy-tailed case}
\label{S:heavy}
In this section, we prove Theorem~\ref{theo: B} and give some properties of the limiting genealogical process, which are of
independent interest. Recall that the fraction of individuals affected by an event is set constant, equal to $u\in (0,1]$, and the radii of the events are sampled according to the intensity measure
$$
\mu(dr)= r^{-\alpha-d-1}\ind_{\{r\geq 1\}}\, dr,
$$
where $d$ is the dimension of the geographical space.

As in the proof of Theorem~\ref{theo: A}, due to the duality relations (\ref{eq duality}) we need only establish
the asymptotic behaviour of the rescaled genealogical process $(\cA^n_t)_{t\geq 0}$ of a finite sample of individuals,
defined in our previous notation by
$$
\cA^n_t\equiv \big(\xi^{n,1}_t,\ldots,\xi^{n,N^n_t}_t\big):= \big(n^{-1/\alpha}\xi^1_{nt},\ldots,n^{-1/\alpha}\xi^{N_{nt}}_{nt}\big).
$$
In words, we speed up time by a factor $n$ and scale down the spatial locations of the lineages by $n^{1/\alpha}$.
Indeed, if we can show that the finite-dimensional distributions of $\cA^n$ converge to those of a system of coalescing
processes $\cA^{\infty}$ that has sufficiently nice properties (i.e., which can be used to construct a dual
$\Xi$-valued process $\rho^{(\alpha)}$ using the technique of \cite{EVA1997}), then the same arguments as those used in the proof of
Theorem~\ref{theo: A} will grant us the convergence of the finite-dimensional distributions of $\rho^n$ to those of $\rho^{(\alpha)}$.
Then it will remain to show that $\cA^{\infty}$ satisfies the conditions of Lemma~\ref{lem: density}$(i)$ to obtain the desired form
for the local densities of $1$'s, $\wa(t,x)$, and to use $(\ref{point duality})$ to characterize the correlations between these
Bernoulli random variables.
Hence, the crucial step is to prove the following proposition.
\begin{prop}\label{prop: geneal B}
There exists a system $\cA^{\infty}$ of coalescing symmetric $\alpha$-stable L\'evy processes such that
$$
\cA^n\rightarrow \cA^{\infty},\qquad \mathrm{as}\ n\rightarrow \infty,
$$
in the sense of weak convergence of the finite-dimensional distributions. Moreover, if we define the process
$\hat{\cA}^n$ and $\hat{\cA}^{\infty}$ in an analogous way to the corresponding processes
in Lemmas~\ref{lem: geneal A1} and \ref{lem: geneal A2}, we also have convergence of the one-dimensional
distributions of $\hat{\cA}^n$ to those of $\hat{\cA}^{\infty}$.
\end{prop}

\medskip
\noindent{\bf Proof of Proposition~\ref{prop: geneal B}.} Our aim is to write down the generator $\cG^n$ of $\cA^n$, and to show that
it converges to the generator of a system of coalescing symmetric $\alpha$-stable processes. Up to now, we were able to be
rather vague about the precise representation of the ancestral lineages, but in order to write down a sensible generator
we now need to be more precise. Suppose we start with $k$ lineages. The system at any time $t\geq 0$ is represented by a
marked partition of $\{1,\ldots,k\}$.  Each block of $\cA^n_t$ contains the labels of all individuals in the initial sample which
have the same ancestor at time $t$ in the past (that is, whose ancestral lineages merged before $t$), and the mark associated to the
block gives the spatial location of this ancestor at time $t$.

Since only the lineages present in the area hit by an event can be affected by this event, for every $y\in \R^d$, $r>0$ and
every marked partition $A$ let us write $J(y,r,A)$ for the set of indices of lineages (blocks) of $A$ whose mark belongs to
$B(y,r)$ (to index the blocks of $A$, we rank them in increasing order of the smallest label that each contains).
For convenience, we shall also use the notation $J_n(y,r,A):=J(n^{-1/\alpha}y,n^{-1/\alpha}r,A)$. Next, if $A$ contains $m$
blocks and $I\subset \{1,\ldots,m\}$, then for every $z\in \R^d$ we write $\Phi_I(A,z)$ for the marked partition obtained by
merging all blocks of $A$ indexed by $i\in I$ and by assigning the mark $z$ to this new block
(the other blocks and marks remain unchanged).
For instance, if $A=\{(\{1,5\},x_1),(\{2,3\},x_2),(\{4,6\},x_3),(\{7\},x_4)\}$ and $I=\{1,4\}$, then
$$
\Phi_I(A,z)= \big\{(\{1,5,7\},z),(\{2,3\},x_2),(\{4,6\},x_3)\big\}.
$$
Finally, we write $|I|$ for the cardinality of the set $I$, and we recall that $V_r$ denotes the volume of a ball of radius $r$.

Because lineages jump and merge at finite rate,
the generator $\cG$ of the system of unrescaled lineages $(\cA_t)_{t\geq 0}$ can be expressed as follows.  For every bounded
measurable function $f$ and every marked partition $A$ (of some finite set $\{1,\ldots,k\}$),
\begin{equation}\label{generator A}
\cG f(A)=\int_{\R^d}dy\int_0^{\infty}\mu(dr)\int_{B(y,r)}\frac{dz}{V_r}\sum_{I\subset J(y,r,A)}u^{|I|}(1-u)^{|J\setminus I|}\big[f(\Phi_I(A,z))-f(A)\big],
\end{equation}
where in the above and what follows we write $J$ as a shorthand notation for $J(y,r,A)$. Indeed, if an event occurs in $B(y,r)$ and the parent is chosen at location $z$, then every lineage present in this area is
affected by the event with probability $u$, independently of each other, and all lineages that are affected merge and jump
onto the location $z$ of their parent.

Mutiplying time by $n$ and marks by $n^{-1/\alpha}$, we obtain from the expression in (\ref{generator A}) that the
generator of $\cA^n$ is given, for every $f$ and $A$ as above, by
$$
\cG^n f(A)=n \int_{\R^d}dy\int_0^{\infty}\mu(dr)\int_{B(y,r)}\frac{dz}{V_r}\sum_{I\subset J_n(y,r,A)}u^{|I|}(1-u)^{|J\setminus I|}\big[f(\Phi_I(A,n^{-1/\alpha}z))-f(A)\big].
$$
To see where the sum comes from, observe that an unrescaled mark belongs to $B(y,r)$ iff its rescaled version belongs to
$B(n^{-1/\alpha}y,n^{-1/\alpha}r)$, and that the affected (rescaled) lineages jump onto $n^{-1/\alpha}z$ when their unrescaled
counterparts jump to $z$. Making the change of variables $z'=n^{-1/\alpha}z$, and then $y'=n^{-1/\alpha}y$ and $r'=n^{-1/\alpha}r$,
we obtain that $\cG^n(A)$ is equal to \setlength\arraycolsep{0pt}
\begin{align}
n^{1+\frac{d}{\alpha}}& \int_{\R^d}dy\int_1^{\infty}\frac{dr}{r^{\alpha+d+1}}\int_{B(n^{-1/\alpha}y,n^{-1/\alpha}r)}\frac{dz}{V_r}\sum_{I\subset J_n(y,r,A)}u^{|I|}(1-u)^{|J\setminus I|}\big[f(\Phi_I(A,z))-f(A)\big]\nonumber\\
= & \int_{\R^d}dy\int_{n^{-1/\alpha}}^{\infty}\frac{dr}{r^{\alpha +d+1}}\int_{B(y,r)}\frac{dz}{V_r}\sum_{I\subset J(y,r,A)}u^{|I|}(1-u)^{|J\setminus I|}\big[f(\Phi_I(A,z))-f(A)\big]\nonumber\\
= & \int_{\R^d}dy\int_{n^{-1/\alpha}}^{\infty}\frac{dr}{r^{\alpha +d+1}}\int_{B(y,r)}\frac{dz}{V_r}\sum_{I\subset J(y,r,A), |I|\geq 2}u^{|I|}(1-u)^{|J\setminus I|}\big[f(\Phi_I(A,z))-f(A)\big]\nonumber \\
& + \int_{\R^d}dy\int_{n^{-1/\alpha}}^{\infty}\frac{dr}{r^{\alpha +d+1}}\int_{B(y,r)}\frac{dz}{V_r}\sum_{i\in J(y,r,A)}u(1-u)^{|J|-1}\big[f(\Phi_{\{i\}}(A,z))-f(A)\big].\label{split gen}
\end{align}
Let us define $\delta(A)$ as half of the minimal pairwise distance between marks in $A$ ($\delta(A):=+\infty$ if $A$ contains only one block),
and let us show that for every $A$ such that $\delta(A)>0$ and every $f$ compactly supported and of class $C^2$ with respect to the
marks, $\cG^nf(A)$ converges as $n\rightarrow \infty$ towards the quantity $\cGa f(A)$ defined by
\begin{align}
&\cGa f(A) \nonumber\\
&:= \int_{\R^d}dy\int_0^{\infty}\frac{dr}{r^{\alpha +d+1}}\int_{B(y,r)}\frac{dz}{V_r}\sum_{I\subset J(y,r,A), |I|\geq 2}u^{|I|}(1-u)^{|J\setminus I|}\big[f(\Phi_I(A,z))-f(A)\big]\nonumber\\
& \quad + u\sum_{i=1}^{|A|}\int_{\R^d}dy\int_0^{\infty}\frac{\ind_{\{x_i\in B(y,r)\}}dr}{r^{\alpha +d+1}}(1-u)^{|J(y,r,A)|-1}\nonumber \\
& \qquad \qquad \qquad \qquad \qquad \qquad \times \int_{B(y,r)}\frac{dz}{V_r}\big[f(\Phi_{\{i\}}(A,z))-f(A)-\la z-x_i,\nabla_i f(A)\ra\ind_{\{|z-x_i|\leq 1\}}\big] \nonumber\\
& \quad + u\sum_{i=1}^{|A|}\int_{\R^d}dy\int_0^{\infty}\frac{\ind_{\{x_i\in B(y,r)\}}dr}{r^{\alpha+d+1}}(1-u)^{|J(y,r,A)|-1}\int_{B(y,r)}\frac{dz}{V_r}\la z-x_i,\nabla_i f(A)\ra\ind_{\{|z-x_i|\leq 1\}}, \label{limit gen}
\end{align}
where $|A|$ denotes the number of blocks of $A$, $x_i$ is the mark of the $i$-th block, $\nabla_i f$ is the gradient of $f$ with
respect to $x_i$ and $\la \cdot,\cdot\ra$ is the scalar product in $\R^d$.

We shall comment on the different terms of $\cGa f(A)$ later. For now, let us show the desired convergence, as well as the
finiteness of $\cGa f(A)$. Let us start with the first term on the right-hand side of (\ref{split gen}). By definition
of $\delta(A)$, a ball of radius $r<\delta(A)$ cannot contain more than $1$ lineage (mark), so that the integral over $r$ runs in
fact from $n^{-1/\alpha}\vee \delta(A)$ to $+\infty$. For $n$ large enough, this first term is thus equal to
$$
\int_{\R^d}dy\int_{\delta(A)}^{\infty}\frac{dr}{r^{\alpha +d+1}}\int_{B(y,r)}\frac{dz}{V_r}\sum_{I\subset J(y,r,A), |I|\geq 2}u^{|I|}(1-u)^{|J\setminus I|}\big[f(\Phi_I(A,z))-f(A)\big],
$$
and so is the first term of $\cGa f(A)$. Since $u\in (0,1]$, $f$ is bounded, the sum over $I$ is finite and since any event
location $B(y,r)$ must intersect the compact support of $f$ to have a nonzero contribution to the generator (so that we may restrict
the integral over $y$ to some ball $B(0,r+ \Delta(f))$ with $\Delta(f)$ depending only on $f$), there exists a constant $C(f)>0$,
independent of $A$, such that the absolute value of the first term of $\cGa f(A)$ is bounded by
\begin{equation}\label{boundedness}
C(f)\, 2^{|A|}\int_{\delta(A)}^{\infty}\frac{dr}{r^{\alpha +d+1}}\, r^d <\infty.
\end{equation}

Now consider the second term on the right-hand side of~(\ref{split gen}).  Let us split it once again into
\begin{align}
& \int_{\R^d}dy\int_{n^{-1/\alpha}}^{\infty}\frac{dr}{r^{\alpha +d+1}}\int_{B(y,r)}\frac{dz}{V_r}\sum_{i\in J(y,r,A)}u(1-u)^{|J|-1}\nonumber \\
& \qquad \qquad \qquad \qquad \qquad \qquad \times \big[f(\Phi_{\{i\}}(A,z))-f(A)-\la z-x_i,\nabla_i f(A)\ra\ind_{\{|z-x_i|\leq 1\}}\big] \label{G1}\\
& \quad + \int_{\R^d}dy\int_{n^{-1/\alpha}}^{\infty}\frac{dr}{r^{\alpha+d+1}}\int_{B(y,r)}\frac{dz}{V_r}\sum_{i\in J(y,r,A)}u(1-u)^{|J|-1}\la z-x_i,\nabla_i f(A)\ra\ind_{\{|z-x_i|\leq 1\}} \label{G2}.
\end{align}
We rewrite $\sum_{i\in J(y,r,A)}$ as $\sum_{i=1}^{|A|}\ind_{\{x_i \in B(y,r)\}}$, and, for $n$ large enough, we split the integral
over $r\in [n^{-1/\alpha},\infty)$ in~(\ref{G2}) into the integral over $[n^{-1/\alpha},\delta(A))$
and that over $[\delta(A),\infty)$.  The second integral is finite for the  same reasons as in (\ref{boundedness}).
On the other hand, if $r<\delta(A)$ then $J(y,r,A)\leq 1$ for every $y$, and so the first integral is equal to
\begin{align*}
u \sum_{i=1}^{|A|} \int_{\R^d}dy &\int_{n^{-1/\alpha}}^{\delta(A)}\frac{dr}{r^{\alpha+d+1}}\, \ind_{\{x_i\in B(y,r)\}}\int_{B(y,r)}\frac{dz}{V_r}\, \la z-x_i,\nabla_i f(A)\ra \ind_{\{|z-x_i|\leq 1\}}\\
= & \ u\sum_{i=1}^{|A|} \int_{n^{-1/\alpha}}^{\delta(A)}\frac{dr}{V_r r^{\alpha+d+1}}\int_{B(x_i,1)}dz\int_{\R^d}dy\ \ind_{\{|x_i-y|\leq r\}}\ind_{\{|z-y|\leq r\}} \la z-x_i,\nabla_i f(A)\ra \\
= & \ u\sum_{i=1}^{|A|} \int_{n^{-1/\alpha}}^{\delta(A)}\frac{dr}{V_r r^{\alpha+d+1}}\int_{B(x_i,1)}dz\ \Big(\mathrm{Vol}\big(B(x_i,r)\cap B(z,r)\big)\Big) \la z-x_i,\nabla_i f(A)\ra,
\end{align*}
and, by symmetry, the integral over $z$ is equal to $0$ for every $r$. The integral in (\ref{G2}) is thus equal to
$$
u\sum_{i=1}^{|A|}\int_{\R^d}dy\int_{\delta(A)}^{\infty}\frac{\ind_{\{x_i\in B(y,r)\}}dr}{r^{\alpha+d+1}}(1-u)^{|J(y,r,A)|-1}\int_{B(y,r)}\frac{dz}{V_r}\la z-x_i,\nabla_i f(A)\ra\ind_{\{|z-x_i|\leq 1\}},
$$
and if we decompose the range $(0,\infty)$ over which we integrate $r$
in the third term of $\cGa f(A)$ into $(0,\delta(A))$ and $[\delta(A),\infty)$, we find that the integral over the
latter is equal to the quantity above.

Finally, let us show that (\ref{G1}) converges to the second term of $\cGa f(A)$. This time, we split (\ref{G1}) into
\begin{align*}
u &\sum_{i=1}^{|A|}\int_{\R^d} dy \int_{n^{-1/\alpha}}^{\infty}\frac{\ind_{\{x_i\in B(y,r)\}}dr}{r^{\alpha + d+1}}(1-u)^{|J(y,r,A)|-1}\int_{B(y,r)}\frac{dz}{V_r}\big(f(\Phi_{\{i\}}(A,z))-f(A)\big)\ind_{\{|z-x_i|> 1\}} \\
+ & u \sum_{i=1}^{|A|}\int_{\R^d} dy \int_{n^{-1/\alpha}}^{\infty}\frac{\ind_{\{x_i\in B(y,r)\}}dr}{r^{\alpha + d+1}}(1-u)^{|J(y,r,A)|-1}\\
& \qquad \qquad \qquad \qquad \qquad \qquad \times \int_{B(y,r)}\frac{dz}{V_r}\big(f(\Phi_{\{i\}}(A,z))-f(A)- \la z-x_i,\nabla_i f(A)\ra\big)\ind_{\{|z-x_i|\leq  1\}}.
\end{align*}
The first term is finite for the same reasons as in (\ref{boundedness}), since for the parent to be at distance greater
than $1$ from the affected lineage, one must have $r>1/2$. Now, using the same steps as above, we obtain that the
second term is equal to
\begin{align}
& u\sum_{i=1}^{|A|} \int_{B(x_i,1)}dz\int_{n^{-1/\alpha}\vee \frac{|z-x_i|}{2}}^{\infty}\frac{dr}{V_r r^{\alpha + d+1}}\int_{B(z,r)\cap B(x_i,r)}dy\ (1-u)^{|J(y,r,A)|-1}\nonumber \\
& \qquad \qquad\qquad\qquad \qquad\qquad \qquad\qquad\qquad \times \big(f(\Phi_{\{i\}}(A,z))-f(A)- \la z-x_i,\nabla_i f(A)\ra\big).\label{G3}
\end{align}
But $f$ is of class $C^2$ and has compact support, and so we can find a constant $\tilde C(f)>0$, independent of $A$, such that for
every $i$ and every $z\in B(x_i,1)$,
$$
\big|f(\Phi_{\{i\}}(A,z))-f(A)- \la z-x_i,\nabla_i f(A)\ra\big| \leq \tilde C(f)|z-x_i|^2.
$$
As a consequence, the absolute value of the quantity in (\ref{G3}) is bounded by \setlength\arraycolsep{1pt}
\begin{align}
u\tilde C & \sum_{i=1}^{|A|} \int_{B(x_i,1)}dz\int_{n^{-1/\alpha}\vee \frac{|z-x_i|}{2}}^{\infty}\frac{dr}{r^{\alpha + d+1}}\frac{\mathrm{Vol}\big(B(z,r)\cap B(x_i,r)\big)}{V_r}\, |z-x_i|^2 \nonumber\\
&\leq  uC'|A| \int_{B(0,1)} dz\ |z|^2\big(n^{-1/\alpha}\vee (|z|/2)\big)^{-\alpha-d}\nonumber \\
& = uC'|A|\left\{n^{1+\frac{d}{\alpha}}\int_{B(0,2n^{-1/\alpha})}dz\ |z|^2 + 2^{\alpha+d}\int_{B(0,1)\setminus B(0,2n^{-1/\alpha})} dz\ |z|^{2-\alpha-d}\right\} \nonumber\\
& \leq C''|A| \big\{n^{-\frac{2-\alpha}{\alpha}} + C'''\big(1- n^{-\frac{2-\alpha}{\alpha}}\big)\big\},\label{error gen}
\end{align}
where all the constants appearing in this bound depend on $f$, $d$ and $\alpha$, but not on $A$. Since $\alpha<2$, (\ref{G3})
remains bounded as $n\rightarrow \infty$ and (\ref{G1}) converges to
\begin{align*}
u\sum_{i=1}^{|A|}\int_{\R^d}dy\int_0^{\infty}\frac{\ind_{\{x_i\in B(y,r)\}}dr}{r^{\alpha +d+1}}(1-u)^{|J(y,r,A)|-1}\int_{B(y,r)}\frac{dz}{V_r}\big[&f(\Phi_{\{i\}}(A,z))-f(A)\\
& -\la z-x_i,\nabla_i f(A)\ra\ind_{\{|z-x_i|\leq 1\}}\big],
\end{align*}
which is precisely the second term of $\cGa f(A)$ (and is finite according to the analysis above).
Tracing back our calculations, we see that for $n$ large enough (such that $n^{-1/\alpha}<\delta(A)$) the difference between
$\cG^n f(A)$ and $\cGa f(A)$ is equal to the difference between the quantity in (\ref{G3}) and its counterpart in $\cGa f(A)$
(that is, the second term of $\cGa f(A)$ in which $y$ is only integrated over $B(x_i,1)$). Hence, according to~(\ref{error gen}),
for every $n> \delta(A)^{-\alpha}$
$$
\big|\cG^n f(A)-\cGa f(A)\big| \leq c_f\, |A|\,n^{-\frac{2-\alpha}{\alpha}} ,
$$
where the constant $c_f$ is again independent of $A$. Consequently, for every $f$ which is compactly supported and of class $C^2$ with
respect to the marks, the function $\cGa f$ is bounded and the convergence
\begin{equation}\label{conv Gn}
\lim_{n\rightarrow \infty}\ \sup_{\delta(A)>\e,|A|\leq k}\ \big|\cG^n f(A) - \cGa f(A)\big| = 0
\end{equation}
holds for any choice of $\e>0$ and $k\in \N$.

To conclude the proof of Proposition~\ref{prop: geneal B}, let us use the following result, whose proof we postpone for the sake of clarity. For every $\e>0$, let $t_\e$ be the first time at which at least two lineages lie at distance less than $\e>0$ without having coalesced.
\begin{lemma}\label{lem: process exists}
For every initial value $A_0$ such that $\delta(A_0)>0$, we have
\begin{equation}\label{exists limit}
\lim_{\e\rightarrow 0}\P_{A_0}[t_{\e}<\infty] = 0.
\end{equation}
As a consequence, the martingale problem associated to $(\cGa,A_0)$ has a unique solution (with c\`adl\`ag paths) for any initial value $A_0$ satisfying $\delta(A_0)>0$. Let us denote this solution by $\cA^{\infty}$. Then $\cA^{\infty}$ is a consistent system of coalescing symmetric $\alpha$-stable processes.
\end{lemma}

Let us suppose that Lemma~\ref{lem: process exists} has been established, and verify that the conditions of
Theorem~4.8.2$(b)$ of \cite{EK1986} are then fulfilled. First, one can check that the set of functions $f$ considered above is
dense in the set of all bounded continuous functions on marked partitions. We can thus restrict our attention to these particular
functions. Second, (\ref{exists limit}) enables us to use (\ref{conv Gn}) and dominated convergence to obtain that
Condition (8.7) of Theorem~4.8.2$(b)$ of \cite{EK1986} is satisfied, and consequently that the finite-dimensional
distributions of $\cA^n$ converge weakly to those of $\cA^{\infty}$ as $n$ tends to infinity. The arguments for
the convergence of the one-dimensional distributions of $\hat\cA^n$ are the same as in the case with fixed radii,
and so the proof of Proposition~\ref{prop: geneal B} is now complete.
\hfill $\Box$

\medskip
Before proving Lemma~\ref{lem: process exists}, let us study some of properties of the `genealogical' process $\cA^{\infty}$. Indeed, in order to use Lemma~\ref{lem: density}$(i)$, we need to show that (\ref{cond Bernoulli}) holds. In fact we can be more precise about the way coalescence occurs.

\begin{lemma}\label{L:coalH}
Sample two individuals at separation $x$, and consider their ancestral lineages $(X_t,t \ge 0)$, $(Y_t,t\ge 0)$. Let
$$
\tau=\inf\{t\geq 0: X_s = Y_s \text{ for all } s \geq t\}
$$
be their coalescence time. Then $\tau < \infty$ almost surely, and moreover, there exists a random variable $Z$,
a.s. finite and independent of $x$, such that
\begin{equation} \label{Coalscal}
\tau \preceq x^\alpha Z,
\end{equation}
where $\preceq$ stands for stochastic domination.
\end{lemma}

\noindent {\bf Proof of Lemma \ref{L:coalH}.}
In essence, the strategy of the proof consists of showing that if the two lineages start at distance $a>0$,
they have some positive chance (independent of $a$) of coalescing before they either separate to a distance greater than $2a$
or come within distance less than $a/2$ of each other. The dependence on $x^{\alpha}$ in the lemma then comes from the fact
that the time needed to coalesce, or separate, or get closer by a factor of $2$, is of the order of $x^{\alpha}$
when the initial separation is $x$.

By translation invariance, we may assume without loss of generality that the origin of $\R^d$ sits at the midpoint between $X_0$ and $Y_0$. Let $T(x)$ be the first time that \emph{any} point in $B:=B(0,x)$ is touched by an event whose radius $r$ is greater than $x/4$. Then we claim that $T(x)$ is an exponential random variable whose rate $\lambda(x)$ is given for every $x>0$ by
\begin{equation}\label{value lambda}
\lambda(x) = \int_{x/4}^{\infty} \frac{d\ell}{\ell^{d+1+\alpha}}\, \hbox{Vol}(B(0,x+\ell)).
\end{equation}
Indeed, recall the intensity measure (\ref{def mu B}) we introduced before rescaling the process. In
the original units of time and space, the rate at which any point of the closed ball $B(0,x)$ ($x\geq 4$) is hit by an event of radius greater than $x/4$ is given by
$$
\int_{\R^d}dz \int_{x/4}^{\infty} \frac{d\ell}{\ell^{d+1+\alpha}}\, \ind_{\{B(0,x)\cap B(z,\ell)\neq \emptyset\}} = \int_{x/4}^{\infty}\frac{d\ell}{\ell^{d+1+\alpha}}\, \hbox{Vol}(B(0,x+\ell)).
$$\
Multiplying this rate by $n$ and looking at distances of the form $xn^{1/\alpha}$, a simple change of variables gives us that for every $x\geq 4n^{-1/\alpha}$, the rescaled rate of interest is also equal to the expression above, independently of $n$. Passing to the limit $n\rightarrow \infty$ yields (\ref{value lambda}).

Now, setting $\ell = rx$ we can write
\begin{align}
\lambda(x) &= x^{-d-\alpha} \int_{1/4}^{\infty} \frac{dr}{r^{d+1+\alpha}}\, \hbox{Vol}(B(0,x+rx)) \nonumber \\
& = x^{-\alpha}\int_{1/4}^{\infty} \frac{dr}{r^{d+1+\alpha}}\, \hbox{Vol}(B(0,1+r)) = Cx^{-\alpha},\label{RateT(x)}
\end{align}
where the constant $C$ is independent of $x$.

On the other hand, similar calculations enable us to see that the rate at which $B$ is
entirely contained within the area $B(z,r)$ of an event is given by
\begin{align*}
\int_{\R^d}dz \int_{|z|+x}^{\infty} \frac{d\ell}{\ell^{d+1+\alpha}} & = \int_x^{\infty} \frac{d\ell}{\ell^{d+1+\alpha}} \int_{\R^d}dz\, \ind_{\{|z|\leq \ell-x\}} \\
& = x^{-\alpha} \int_1^{\infty} \frac{dr}{r^{d+1+\alpha}}\, \hbox{Vol}(B(0,r-1)) = C'x^{-\alpha},
\end{align*}
where we used the same change of variable as before and $C'>0$ is again independent of $x$. As a consequence, with probability $p_0 := C'/C$ independent of $x$, the first event of radius greater than $x/4$ that hits at least one point of $B$ actually covers the whole ball. Moreover, \eqref{RateT(x)} also implies that for arbitrary $q\geq 1/4$, the radius $R(x)$ of the event occurring at time $T(x)$ satisfies
\begin{equation}\label{R(x)unif_ub}
\P( R(x) > qx) \leq c q^{-\alpha},
\end{equation}
for some constant $c$ which does not depend on $x$ or $q$.

Let $\tilde X, \tilde Y$ be the motion of the lineages as governed by all the events except those that affect some point in $B$ and whose radius is greater than $x/4$. Then by the Poisson point process formulation of the reproduction events, $T(x)$ is independent of $\tilde X, \tilde Y$ and $(X_t, Y_t, t <T(x))$ coincides with $(\tilde X_t, \tilde Y_t, t < T(x))$. Let $S(x) := \inf\{ t\ge 0: \tilde D_t \le x/2 \text{ or }\tilde X_t \notin B \text{ or }\tilde Y_t \notin B\}$, where $D_t = |\tilde X_t - \tilde Y_t|$. Fix $\delta>0$, and define the following events:
$$
E:= \{T(x) \le \delta x^\alpha\}, \qquad \qquad F:= \{S(x) \ge \delta x^\alpha\}.
$$
Then $E$ and $F$ are independent, and by \eqref{RateT(x)} there exists $p(\delta)>0$ such that $\P(E) = p(\delta)$ for all $x>0$.
A similar property holds for $F$. Indeed, note first that up until the time $S(x)$, the trajectories $\tilde X$ and $\tilde Y$ are
independent, since the trajectories can only move as a result of events occurring in necessarily disjoint regions of space.
Moreover, it is easy to check that
\begin{equation}\label{eqXtildeYtilde}
\Big(\frac1x\,\tilde X_{t x^\alpha \wedge S(x) }, \, \frac1x\,\tilde Y_{t x^\alpha \wedge S(x)}\Big)_{t \ge 0}
\end{equation}
has the same distribution as the pair $(\tilde X_{t\wedge S(1)},\tilde Y_{t\wedge S(1)})_{t\ge 0}$ obtained by taking $x=1$:
both coordinates of this process perform independent stable L\'evy processes where each jump greater than $1/4$ occurring
in $B(0,1)$ is removed, and the process is stopped when either coordinate leaves $B(0,1)$ or they come within distance $1/2$ of one
another. (A formal proof is given by comparing the generators; the generator of the pair $(\tilde X, \tilde Y)$ is the same as $\cG^\alpha f(A)$ when $A$ has two blocks, as defined in \eqref{limit gen}, but with the first term equal to 0).

Hence for all $x>0$, $\P(S(x) \ge \delta x^\alpha) = \P(S(1) \ge \delta)=: q(\delta)$, and $q(\delta)>0$ whenever $\delta$
is chosen small enough.

Let us denote the centre, radius and impact parameter of the event taking place at time $T(x)$ by $(Z(x),R(x),u)$. We shall
say that a \emph{success} occurs if both $E$ and $F$ occur, and if
\begin{enumerate}
\item[(a)] $B(0,x) \subset B(Z(x), R(x))$,
\item[(b)] both $X_{T(x)}, Y_{T(x)}$ are affected by the event occurring at time $T(x)$ (this is possible since under these assumptions, $X_{T(x)}$ and $Y_{T(x)}$ are still both in $B(0,x)$ which is entirely covered by the event.)
\end{enumerate}
Note that by the above discussion,
\begin{equation}\label{success}
\wp :=\P(\text{ success }) = p(\delta)q(\delta) p_0 u^2,
\end{equation}
independently of $x>0$.

If a success did not occur, we say that a \emph{failure} has occurred. Since the success probability is independent of $x$ and
the waiting time between two attempts is always stochastically bounded by an exponential random variable of the form
$T(y)$ (which is a.s.~finite), we deduce that after a Geometric$(\wp)$ number $N$ of attempts, success is guaranteed,
hence $\tau< \infty$ almost surely. Moreover,
in the case of failure, consider the mutual distance $D_{T(x) \wedge S(x)}$ between the two lineages at time $T(x) \wedge S(x)$.
Then $D_{T(x) \wedge S(x)} \le 2x+R(x)$. From \eqref{R(x)unif_ub} we can deduce that there exists a random variable $R$,
independent of $x$ and a.s. finite, such that $2+R(x)/x \preceq R$ in the sense of stochastic domination.
Let $R_1, R_2, \ldots$ be a sequence of i.i.d. random variables with distribution $R$.
The strong Markov property and (\ref{RateT(x)}) then show that
$$
\tau \preceq x^{\alpha}\Big\{\mathcal{E}[C]+ \mathcal{E}\big[C R_1^{-\alpha}\big] + \ldots + \mathcal{E}\big[C (R_1 \cdots R_N)^{-\alpha}\big]\Big\},
$$
where $\mathcal{E}[y]$ stands for an exponential random variable with parameter $y$ and all the above exponential
random variables are conditionally independent given their arguments. Define $Z$ as the random variable within the
curly brackets to conclude. \hfill $\Box$

\begin{rmk}
\label{rk: MRCA}
The system $\cA^{\infty}$ inherits the consistency property from its construction as the limit of $\cA^n$ (this property can also be
shown directly from the generator of $\cA^{\infty}$). Hence, a notable consequence of Lemma~\ref{L:coalH} is that any finite sample of
lineages finds its most recent common ancestor in finite time with probability one. The same kind of behaviour, as well as the
convergence of the forwards-in-time process to a field of correlated Bernoulli random variables, was already observed by Evans
in the case where the genealogical process of his \emph{continuous sites stepping-stone model} is a system of one-dimensional
independent $\alpha$-stable motions coalescing instantly upon meeting. See \S5 in \cite{EVA1997} for a full description of his results.
However, the underlying mechanisms are quite different here. Not only
does Lemma~\ref{L:coalH} hold for any $\alpha\in (1,2)$ and any dimension, which cannot be the case in Evans' framework since two
independent stable processes may not meet, but even in dimension $1$ the way lineages coalesce is different:
the limit in (\ref{exists limit}) shows that two lineages of $\cA^{\infty}$ have no chance to meet, but their coalescence is due to
the fact that large events of the appropriate size are just frequent enough to catch them even when they are very far from each
other. As a last consequence, it is then possible to see multiple mergers during the evolution of $\cA^{\infty}$, which is not the
case when the $\alpha$-stable processes move independently of each other and coalesce only when they meet.
\end{rmk}

Let us now finish with the proof of Lemma~\ref{lem: process exists} and of Theorem~\ref{theo: B}.  Recall that for any
marked partition $A$, $\delta(A)$ stands for half the minimum distance between two marks in $A$ ($\delta(A)=+\infty$ if $A$ has only one block).

\medskip
\noindent{\bf Proof of Lemma \ref{lem: process exists}.}
Because most of the ideas and computations we shall use to establish (\ref{exists limit}) are developed in detail in the proof of
Lemma~\ref{L:coalH}, we only present an outline here and refer to that proof for more precise arguments.
Since we always deal with partitions of some finite set, it is sufficient to show the result when $A_0$ consists of just
two blocks starting at some positive separation.

If $x>0$ denotes the initial distance between our two lineages, let us call $T(x)$ the first time at which any of the lineages is in
the geographical area of an event of radius greater than $x/4$, and let us call $S(x)$ the first time at which the distance between
the two lineages is greater than $2x$, or less than $x/2$. Notice that the lineages evolve independently until the random time
$T(x)\wedge S(x)$, since they are hit by events that are necessarily disjoint until that time. Moreover, they both move according to
the law of a symmetric $\alpha$-stable process whose large jumps have been truncated (see reasoning below \eqref{eqXtildeYtilde}). Hence it is not
difficult to show that $S(x)$ is of the order of $x^{\alpha}$, and so is $T(x)$, while the coalescence rate of two lineages at
distance $x$ is commensurate with $x^{-\alpha}$. Using the more careful analysis performed in the proof of Lemma~\ref{L:coalH},
we can in fact conclude that the probability $p_0$ that the two lineages coalesce before their distance doubles or is divided by two is
not only positive, but also independent of $x$. Together with the fact that $T(y)\wedge S(y)$ is a.s. finite for every $y>0$
(for reasons expounded in Lemma~\ref{L:coalH}), the number of attempts before succeeding to coalesce is a geometric
random variable with parameter $p_0$, which we shall denote by $N$.

As a second step, suppose that the lineages fail to coalesce at time $T(x)\wedge S(x)$. The new location of the
lineage which jumps at that time (at most one of them jumps, otherwise they would coalesce) is uniformly distributed over the area of the event,
and since the lineages are at distance at least $x/2$ from each other just before $T(x)\wedge S(x)$ a small calculation using the
scaling properties of the evolution mechanism shows that the probability $\pi(\eta)$ that their new distance at that time is less
than $\eta x$ satisfies
\begin{itemize}
\item[(a)] $\pi(\eta)$ is independent of $x$,
\item[(b)] $\lim_{\eta\rightarrow 0}\pi(\eta)=0$.
\end{itemize}
As a consequence, if $\eta\in (0,1/10)$ and $k\in \N$, we can write
\begin{equation}\label{lb coal}
\P_{A_0}\big[\hbox{coal. before distance decreases by }\eta^k\big] \geq \E\big[(1-\pi(\eta))^{N-1}\ind_{\{N<k\}}\big].
\end{equation}
Note in passing that, by monotonicity, the same inequality holds if we replace $\eta^k$ by any $\e\leq \eta^k$.

Let us now draw some conclusions from these observations.  We fix $c>0$, and choose $k(c)$ and $\eta(c)$ such that
for every $k\geq k(c)$ and $\eta\leq \eta(c)$,
$$
\P[N\geq k]\leq \frac{c}{2} \qquad \hbox{and}\qquad \E\big[(1-\pi(\eta))^{N-1}\big]\geq 1-\frac{c}{2}.
$$
Then, using the fact that the event described in the left-hand side of (\ref{lb coal}) implies $t_\e =+\infty$ for
every $\e\leq \eta^k x$, we have that, for every such $\e$,
$$
\P_{A_0}[t_\e = \infty]\geq 1-\frac{c}{2} - \frac{c}{2} = 1-c.
$$
Since $c$ was arbitrary, (\ref{exists limit}) follows.

As regards the second part of Lemma~\ref{lem: process exists}, recall from (\ref{limit gen}) that the operator $\cGa$ is defined,
for every function $f$ of class $C^2$ with compact support and every marked partition $A$ satisfying $\delta(A)>0$, by
\begin{align*}
\cGa f(A)=& \int_{\R^d}dy\int_{\delta(A)}^{\infty}\frac{dr}{r^{\alpha +d+1}}\int_{B(y,r)}\frac{dz}{V_r}\sum_{I\subset J(y,r,A), |I|\geq 2}u^{|I|}(1-u)^{|J\setminus I|}\big[f(\Phi_I(A,z))-f(A)\big]\\
& + u\sum_{i=1}^{|A|}\int_{\R^d}dy\int_0^{\infty}\frac{\ind_{\{x_i\in B(y,r)\}}dr}{r^{\alpha +d+1}}(1-u)^{|J(y,r,A)|-1} \\
&  \qquad \qquad \qquad \qquad \times \int_{B(y,r)}\frac{dz}{V_r}\big[f(\Phi_{\{i\}}(A,z))-f(A)-\la z-x_i,\nabla_i f(A)\ra\ind_{\{|z-x_i|\leq 1\}}\big] \\
 + & u\sum_{i=1}^{|A|}\int_{\R^d}dy\int_0^{\infty}\frac{\ind_{\{x_i\in B(y,r)\}}dr}{r^{\alpha+d+1}}(1-u)^{|J(y,r,A)|-1}\int_{B(y,r)}\frac{dz}{V_r}\la z-x_i,\nabla_i f(A)\ra\ind_{\{|z-x_i|\leq 1\}}.
\end{align*}
In particular, if $A=\{(b_1,x_1)\}$ contains only one block and if $f$ is a function of its mark only, then $\cGa f(A)$ is equal to
\begin{align*}
 u&\int_{\R^d}dy\int_0^{\infty}\frac{\ind_{\{x_1\in B(y,r)\}}dr}{r^{\alpha +d+1}} \int_{B(y,r)}\frac{dz}{V_r}\big[f(z)-f(x_1)-\la z-x_1,\nabla f(x_1)\ra\ind_{\{|z-x_1|\leq 1\}}\big] \\
& \qquad \qquad + u \int_{\R^d}dy\int_0^{\infty}\frac{\ind_{\{x_1\in B(y,r)\}}dr}{r^{\alpha+d+1}}\int_{B(y,r)}\frac{dz}{V_r}\la z-x_1,\nabla f(x_1)\ra\ind_{\{|z-x_1|\leq 1\}}\\
= & \ u\int_{\R^d}dz\left(\int_0^{\infty}\frac{dr}{r^{\alpha +d+1}}\frac{\mathrm{Vol}\big(B(z,r)\cap B(x_1,r)\big)}{V_r}\right)\big[f(z)-f(x_1)-\la z-x_1,\nabla f(x_1)\ra\ind_{\{|z-x_1|\leq 1\}}\big] \\
& \qquad \qquad + u \int_{\R^d}dz \left(\int_0^{\infty}\frac{dr}{r^{\alpha +d+1}}\frac{\mathrm{Vol}\big(B(z,r)\cap B(x_1,r)\big)}{V_r}\right) \la z-x_1,\nabla f(x_1)\ra\ind_{\{|z-x_1|\leq 1\}}.
\end{align*}
The second term above is zero (by symmetry) and hence the generator can be written as
\begin{equation}\label{1motion}
\cGa f(A)= \int_{\R^d}dz\ \iota(z - x_1)\ \big(f(z)-f(x_1)\big),
\end{equation}
where the intensity $\iota(w)$ is given by
$$
\iota(w): = u \int_0^{\infty}\frac{dr}{r^{\alpha +d+1}}\frac{\mathrm{Vol}\big(B(w,r)\cap B(0,r)\big)}{V_r}.
$$
Now, one can check that for any $k>0$
$$
k\, \iota(zk^{-1/\alpha})\, d(zk^{-1/\alpha})= \iota(z)\,dz
$$
and so the motion of a single lineage is a symmetric $\alpha$-stable L\'evy process.

When there are at least two blocks, as long as $\delta(\cA^{\infty}_t)>0$ the first term of $\cGa f(\cA^{\infty}_t)$ is
finite and clearly represents the merger and jump at finite rate of several blocks of $\cA^{\infty}$.
However, the coalescence rate of two lineages at distance $\e$ is equal to
$$
u^2\int_{\R^d}dy\int_{\e/2}^{\infty}\frac{dr}{r^{\alpha + d+1}}\ind_{\{x_1,x_2\in B(y,r)\}} = u^2\int_{\e/2}^{\infty}\frac{dr}{r^{\alpha + d+1}}\, \mathrm{Vol}\big(B(x_1,r)\cap B(x_2,r)\big) \propto \e^{-\alpha}
$$
as $\e\rightarrow 0$, and so one can prove the existence of the process $\cA^{\infty}$ only up to $t_\e$, for any $\e>0$.
Yet (\ref{exists limit}) is actually more than what is required to invoke Theorem 4.6.3 in \cite{EK1986} and complete the proof of existence of $\cA^{\infty}$. \hfill $\Box$

\medskip
\noindent {\bf Proof of Theorem \ref{theo: B}.} There is nothing else to do. Duality and the convergence of $\cA^n$ give us the convergence of $\rho^n$ exactly as in the proof of Theorem~\ref{theo: A}. Lemma~\ref{L:coalH} is sufficient to show that (\ref{cond Bernoulli}) holds and so the limiting densities $\wa(t,x)$ are Bernoulli random variables as stated.
\hfill $\Box$

\end{document}